\documentclass[twoside, a4paper, leqno]{article}%{elsarticle}

\usepackage{amsmath, amssymb, amsthm}
\usepackage{mathtools}
\usepackage{mathrsfs}
\usepackage[sort&compress, square, comma, numbers]{natbib} %sort&compress, square, comma, numbers
\usepackage{hyperref}
\usepackage{bm}
\usepackage{graphicx}        % standard LaTeX graphics tool
\usepackage{multicol}        % used for the two-column index
\usepackage[bottom]{footmisc}% places footnotes at page bottom
\usepackage{comment}
\usepackage{geometry}

\graphicspath{   % directories to search the graphics files
{FIGURES/}
}

\begin{document}

\title{Application of the  Finite Element
  Method in a Quantitative Imaging technique}

\author{L. Beilina  \thanks{
Department of Mathematical Sciences, Chalmers University of Technology and
Gothenburg University, SE-42196 Gothenburg, Sweden, e-mail: 
\texttt{larisa@chalmers.se}}}
\date{}

\maketitle

\begin{abstract}
We present the Finite Element Method (FEM) for the numerical solution of the
multidimensional coefficient inverse problem (MCIP) in two
dimensions. This method is used for explicit reconstruction of the
coefficient in the hyperbolic equation using data resulted from a
single measurement. To solve our MCIP we use approximate globally
convergent method and then apply FEM for the resulted equation. Our
numerical examples show quantitative reconstruction of the sound speed
in small tumor-like inclusions.
\end{abstract}

\section{Introduction}

\label{sec:1}

In this work we present the Finite Element Method (FEM) applied for explicit
reconstruction of the coefficient in the hyperbolic equation using
data resulted from a single measurement. This means that the data are
generated by either a single location of the point source or by a
single direction of the incident plane wave. Such multidimensional
coefficient inverse problems (MCIPs) are non-overdetermined ones and
have a lot of applications, such as, e.g., many aspects of acoustics,
electromagnetics, optics, medical imaging, geophysics, etc..

To  solve our MCIP we use approximate globally
convergent method of \cite{BK} where for the solution of MCIP was used
underlying PDE operator instead of least squares functionals. It is
well known that CIPs are both nonlinear and ill-posed.  A main idea of
an approximate globally convergent method is that the least squares
objective functionals are not used in it and the phenomenon of local
minima is avoided.  This method was further verified on
computationally simulated and on experimental data in \cite{BK2, BTKF,
  BTKJ, NBKF} and references therein.

In the current work we apply the finite element method inside
approximate globally convergent method of \cite{BK}. Our goal is obtain
quantitative medical imaging of small inclusions representing
cancerous tumors. This means that we are interested not only in shape
reconstruction but also in the accurate
reconstruction of the contrast of tumor-like inclusions.  Examples of
MCIPs with applications in medicine are inverse problems of magnetic
resonance elastography (MRE) which are studied recently in
\cite{BB, PO2} and references therein. The main feature of this
medical imaging technique is that it allows measure field internally
and this is the case of our numerical examples of section \ref{sec:2}.
We note that for detection of cancer tumors in human
tissue using MRE technique stiffness contrast can be of the order of
2000\% while the density varies only of the order of 8\%
\cite{PO2}. This is the main reason why stiffness is diagnostically
more useful and density is often not considered.

The current work is devoted to the reconstruction of the wave speed in
the wave equation from internal measurements. We consider the
simplified model problem described by the acoustic wave equation
instead of the elastic one.  Application of the method of this work
for another MCIPs can be considered as a topic for a future research.
Numerical examples of section \ref{sec:2} show very accurate and
quantitative reconstruction of tumor-like inclusions which can be even
of the very small sizes (point-size inclusions).  In our future work
we plan to extend the iterative procedure described in this work to
the case of MCIPs with boundary measurements. Similarly with \cite{AB}
an adaptive finite element method can be also considered as a topic for a
future research.

\section{Statements of Forward and Inverse Problems}

\label{sec:2.1}
We consider the Cauchy problem for the hyperbolic equation
\begin{equation}
a(x) u_{tt}=\Delta u\text{ in }\mathbb{R}^{3}\times \left(
0,\infty \right) ,  \label{2.1}
\end{equation}
\begin{equation}
u\left( x,0\right) =0,\text{ }u_{t}\left( x,0\right) =\delta \left(
x-x_{0}\right),  \label{2.2}
\end{equation}
where $\delta$ is the Dirac delta function.  Equation (\ref{2.1})
governs a wide range of applications, including, e.g.  propagation of
acoustic, elastic and electromagnetic waves. In the acoustical case
$c(x)=1/\sqrt{a(x)}$ is the sound speed. In the electromagnetic waves
propagation in a non-magnetic medium, the dimensionless coefficient is
$a(x)=\varepsilon _{r}(x),$ where $\varepsilon _{r}(x)$ is the
spatially distributed dielectric constant of the medium. In the case
of application of equation (\ref{2.1}) in scanning acoustic microscopy
in medical imaging, the sound speed is defined as
$c(x)=\sqrt{(\lambda(x) + 2 \mu(x))/\rho(x)}$, where $\rho(x)$ is the
density and $\lambda(x), \mu(x)$ are the Lam\'e constants of linear
elasticity \cite{BClason}. In the current paper we consider this kind
of applications when the function $a(x)$ in (\ref{2.1}) can be
determined as $a(x)=\rho(x)/(\lambda(x) + 2 \mu(x))$. Then by the
reconstructed function $a(x)$ it will be possible determine the  stiffness
coefficient $\mu(x)$ for the known functions $\rho(x), \lambda(x)$.

Let $\Omega \subset \mathbb{R}^{3}$ be a convex bounded domain with the
boundary $\partial \Omega \in C^{3}.$ Let $d=const.>1.$ We assume that the
coefficient $a(x) $ of equation (\ref{2.1}) is such that
\begin{eqnarray}
a(x) &\in &[1,d],~~a(x) =1\text{ for }x\in \mathbb{R}
^{3}\diagdown \Omega ,  \label{2.3} \\
a(x) &\in &C^{3}\left( \mathbb{R}^{3}\right),  \label{2.4}
\end{eqnarray}
where $d=const.>1$ is \emph{a priori}  known constant.

\textbf{Coefficient Inverse Problem (CIP).} \emph{Suppose that the
  coefficient }$a(x) $\emph{\ satisfies (\ref{2.3}) and
  (\ref{2.4}). Assume that the function }$a(x)$\emph{\ is unknown in
  the domain} $\Omega $\emph{. Determine the function }$a(x)$
\emph{\ for }$x\in \Omega,$\emph{\ assuming that the following
  function }$g\left( x,t\right) $\emph{\ is known for a single source
point  position }$x_{0}\notin \overline{\Omega }$
\begin{equation}
u\left( x,t\right) =g\left( x,t\right)  ~~\forall \left( x,t\right) \in
\partial \Omega \times \left( 0,\infty \right).  \label{2.5}
\end{equation}

In our applications we assume that the source point $x_{0}\notin
\overline{\Omega } $ since we do not want to deal with singularities
near the source location. In real applications the assumption $a(x)
=1$ for $x\in \mathbb{R}^{3}\diagdown \Omega $ means that the function
$a(x) $ has a known constant value outside of the medium of interest
$\Omega.$ The function $g\left( x,t\right) $ in (\ref{2.5}) models
time dependent measurements of the wave field $u(x,t)$ at the boundary of the
domain of interest.

\section{The Transformation Procedure for the Hyperbolic Case}

\label{sec:2.3}

In this section we show how to reduce our Inverse Problem (CIP) to the
Dirichlet boundary value problem for a nonlinear integro-differential
equation. First, we take the Laplace transform of the functions $u$ in the
hyperbolic equation (\ref{2.1}) to get
\begin{equation}
w(x,s)=\int\limits_{0}^{\infty }u(x,t) \mathrm{e}^{-st}dt ~\text{ for }s>\underline{s}%
=const.>0,  \label{2.6}
\end{equation}%
where $\underline{s}$ is a certain number, which we choose in experiments.  It is sufficient to choose $%
\underline{s}$ such that the integral (\ref{2.6}) would converge together
with corresponding $\left( x,t\right) $ derivatives. Thus, we can assume that
the number $\underline{s}$ is sufficiently large. The parameter $s$
is called \emph{pseudo frequency}. It
follows from (\ref{2.1}), (\ref{2.2}), and (\ref{2.6}) that the function $w$
is the solution of the following problem
\begin{equation}
\Delta w-s^{2}a(x) w=-\delta \left( x-x_{0}\right) ,\text{ }x\in
\mathbb{R}^{3},  \label{2.7}
\end{equation}
\begin{equation}
\lim_{\left\vert x\right\vert \rightarrow \infty }w\left( x,s\right) =0,
\label{2.8}
\end{equation}
where the limit in (\ref{2.8}) is proven in \cite{BK}.

We now work only with the function $w(x,s).$ In Theorem 2.7.2 of
\cite{BK} was shown that $w(x,s)>0.$ Hence, we can consider functions
$v(x,s)$  defined as
\begin{equation}\label{transform}
v\left( x,s\right) =\frac{\ln w\left( x,s\right) }{s^{2}}.
\end{equation}

Assuming that the asymptotic behavior in Lemma 2.3 of \cite{BK} holds
we get the following asymptotic behavior of the function $v$
\begin{equation}
\left\Vert D_x^{\beta} D_{s}^k v\left( x,s\right) \right\Vert _{C^{3}\left( \overline{%
\Omega }\right) }=O\left( \frac{1}{s^{k+1}}\right) ,~s\rightarrow \infty
,~~k=0,1.  \label{new3.1}
\end{equation}
Substituting $w=\mathrm{e}^{v}$ in (\ref{2.7}) and noting that the source point
$x_{0}\notin \overline{\Omega }$ and then dividing the resulting
equation for $v$ by $s^{2},$ we obtain
\begin{equation}
\Delta v+s^{2}\left( \nabla v\right) ^{2}=a(x),~~x\in \Omega.  \label{new3.2}
\end{equation}
Denote
\begin{equation}
q\left( x,s\right) =\partial _{s}v\left( x,s\right).  \label{3.6}
\end{equation}%
By (\ref{new3.1}) and (\ref{3.6}) we obtain
\begin{equation*}
v\left( x,s\right) =-\int\limits_{s}^{\infty }q\left( x,\tau \right) d\tau.
\end{equation*}
We rewrite this integral as
\begin{equation}
v\left( x,s\right) =-\int\limits_{s}^{\overline{s}}q\left( x,\tau \right)
d\tau +V\left( x,\overline{s}\right) ,  \label{3.7}
\end{equation}%
where the truncation pseudo frequency $\overline{s}>\underline{s}$ is
a large number. It is important that $V(x,\overline{s})$ in
(\ref{3.7}) is not an arbitrary function, but is defined as
\begin{equation}
V\left( x,\overline{s}\right) =v\left( x,\overline{s}\right) =\frac{\ln
w\left( x,\overline{s}\right) }{\overline{s}^{2}},  \label{2.A}
\end{equation}
where $w\left( x,\overline{s}\right) $ is the Laplace transform (\ref{2.6})
of the solution of the forward problem (\ref{2.1}), (\ref{2.2}) at $s:=
\overline{s}$. The number $\overline{s
}$ should be chosen in numerical experiments. We call the function  $V\left( x,\overline{s}\right)$ as the ``tail'' function  and this function is
unknown. By (\ref{new3.1}) and (\ref{2.A}) we have that
\begin{equation}
\left\Vert V\left( x,\overline{s}\right) \right\Vert _{C^{3}\left( \overline{%
\Omega }\right) }=O\left( \frac{1}{\overline{s}}\right) ,\text{ }\left\Vert
\partial _{\overline{s}}V\left( x,\overline{s}\right) \right\Vert
_{C^{3}\left( \overline{\Omega }\right) }=O\left( \frac{1}{\overline{s}^{2}}%
\right) .  \label{2.B}
\end{equation}
From above equations follows that the tail is small for large values
of $\overline{s}.$ Therefore, one can set $V\left(
x,\overline{s}\right) :=0$. In our recent works \cite{BTKF, BTKJ} we describe
alternative approach  how this tail function can be approximated in
computations.

We now note that in the equation (\ref{new3.2}) the function $a(x)$
does not depends on the parameter $s$.  Thus, differentiating this
equation with respect to $s$ and using (\ref{3.6}) and (\ref{3.7}), we
obtain the following  nonlinear integro-differential equation
\begin{equation}
\begin{split}
& \Delta q-2s^{2}\nabla q\int\limits_{s}^{\overline{s}}\nabla q\left( x,\tau
\right) d\tau +2s\left[ \int\limits_{s}^{\overline{s}}\nabla q\left( x,\tau
\right) d\tau \right] ^{2} \\
& +2s^{2}\nabla q\nabla V-4s\nabla V\int\limits_{s}^{\overline{s}}\nabla
q\left( x,\tau \right) d\tau +2s\left( \nabla V\right) ^{2}=0,~~x\in \Omega .
\end{split}
\label{3.8}
\end{equation}
Conditions(\ref{2.5}) and (\ref{3.6}) imply that we can set the following Dirichlet
boundary condition for the function $q$
\begin{equation}
q\left( x,s\right) =\psi \left( x,s\right) ~~~\forall \left(
x,s\right) \in \partial \Omega \times \left[ \underline{s},\overline{s}%
\right] ,  \label{3.9}
\end{equation}%
where
\begin{equation*}
\psi \left( x,s\right) =\frac{\partial _{s}\ln \varphi }{s^{2}}-\frac{2\ln
\varphi }{s^{3}}
\end{equation*}
and $\varphi \left( x,s\right) $ is the Laplace transform (\ref{2.6}) of the
function $g\left( x,t\right) $ in (\ref{2.5}).

Assume now that we can solve (\ref{3.8})  and find approximations for functions $q$ and $V$ 
in $\Omega $ together with their derivatives $D_{x}^{\alpha
}q,D_{x}^{\alpha }V,\left\vert \alpha \right\vert \leq 2.$ Then the
the  function $a(x)$ can be found
via explicit formula
\begin{equation}
a(x) =\Delta v+\underline{s}^{2}\left( \nabla v\right) ^{2},~~x\in
\Omega,  \label{3.10}
\end{equation}
where the function $v$ can be  obtained via (\ref{3.7}).

\section{The Layer Stripping Procedure}

\label{sec:2.5}

In this section we describe the layer stripping procedure for the
solution of the integro-differential equation (\ref{3.8}).
To do that we make partition of the pseudo frequency interval $[
  \underline{s}, \bar{s}]$ into $N$ sub-intervals
$\bar{s}=s_{0}>s_{1}>\cdots >s_{N}=\underline{s}$ such that
\begin{equation*}
\underline{s}=s_{N}<s_{N-1}<...<s_{1}<s_{0}=\overline{s},s_{i-1}-s_{i}=h,
\end{equation*}
 where $h$ is the step size of every interval and $q\left( x,s\right)
=q_{n}\left( x\right) $ for $s\in (s_{n},s_{n-1}].$ 
Thus, we approximate the function $q\left(
x,s\right) $ in (\ref{3.8}) by a piecewise constant function with
respect to the pseudo frequency $s$. 
We also set
\begin{equation}
q_{0}\equiv 0.  \label{2.99}
\end{equation}%
Hence, integrals in (\ref{3.8})  can be approximated as
\begin{equation}
\int_{s}^{\overline{s}}\nabla q(x,\tau )d\tau =(s_{n-1}-s)\nabla
q_{n}(x)+h\sum_{j=0}^{n-1}\nabla q_{j}(x),s\in (s_{n},s_{n-1}).  \label{4.1}
\end{equation}%
We approximate the boundary condition (\ref{3.9}) by a piecewise constant
function,
\begin{equation}
q_{n}\left( x\right) =\frac{1}{h}\int
\limits_{s_{n}}^{s_{n-1}}\psi \left( x,s\right) ds.  \label{2.37}
\end{equation}
For  every subinterval $\left( s_{n},s_{n-1}\right] ,n\geq 1$ we assume that
functions $q_{j}\left( x\right) ,~j=1,...,n-1, $ for all previous
subintervals are  computed. Then we obtain from (\ref{3.8}) \ the following
  system of approximate equations   for the functions $q_{n}\left( x\right) $
\begin{equation}
\begin{split}
\widetilde{L}_{n}\left( q_{n}\right) & :=\Delta q_{n}-2\left( s^{2}-2s\left(
s_{n-1}-s\right) \right) \left( h\sum\limits_{j=1}^{n-1}\nabla q_{j}\right)
\nabla q_{n} \\
& +2\left( s^{2}-2s\left( s_{n-1}-s\right) \right) \nabla q_{n}\nabla V \\
& =2\left( s_{n-1}-s\right) \left[ s^{2}-s\left( s_{n-1}-s\right) \right]
\left( \nabla q_{n}\right) ^{2}-2sh^{2}\left( \sum\limits_{j=1}^{n-1}\nabla
q_{j}\right) ^{2} \\
& +4s\nabla V\left( h\sum\limits_{j=1}^{n-1}\nabla q_{j}\right)
-2s\left\vert \nabla V\right\vert ^{2},s\in \left( s_{n-1},s_{n}\right] .
\end{split}
\label{2.38}
\end{equation}
The equation (\ref{2.38}) is nonlinear  and this equation depends on the
parameter $s$.  To involve better stability of the computational process, we add the term $-\varepsilon q_{n}$
to the left hand side of equation (\ref{2.38}). Here, $\varepsilon >0$
is a small parameter.
Then we multiply  (\ref{2.38}) by the Carleman Weight Function (CWF)
of the form
\begin{equation}
\mathcal{C}_{n,\lambda }(s)= \mathrm{e}^{\lambda (s-s_{n-1})} ,~~s\in
(s_{n},s_{n-1}],  \label{2.41}
\end{equation}%
and integrate with respect to $s$ over every pseudo frequency interval$(s_{n},s_{n-1}).$ In (\ref{2.41}) the parameter $\lambda \gg 1$ and it should be chosen in numerical
experiments.
Finally, we obtain
\begin{equation}
\begin{split}
L_{n}\left( q_{n}\right) & :=\Delta q_{n}-A_{1,n}\left(
h\sum\limits_{j=0}^{n-1}\nabla q_{j}\right) \nabla q_{n}+A_{1n}\nabla
q_{n}\nabla V-\varepsilon q_{n} \\
& =2\frac{I_{1,n}}{I_{0}}\left( \nabla q_{n}\right) ^{2}-A_{2,n}h^{2}\left(
\sum\limits_{j=0}^{n-1}\nabla q_{j}\left( x\right) \right) ^{2} \\
& +2A_{2,n}\nabla V\left( h\sum\limits_{j=0}^{n-1}\nabla q_{j}\right)
-A_{2,n}\left( \nabla V\right) ^{2},n=1,...,N, 
\end{split}
\label{4.5}
\end{equation}%
with the discretized boundary condition
\begin{equation}
q_{n}(x)=\psi _{n}(x):=\frac{1}{h}\int
\limits_{s_{n}}^{s_{n-1}}\psi (x,\,s)\,\mathrm{d}s\approx \frac{1}{2}[\psi (x,\,s_{n})+\psi (x,\,s_{n-1})],\quad x \in
\partial \Omega .  \label{eq:q7}
\end{equation}
In (\ref{4.5})  coefficients can be computed analytically:
\begin{equation*}
I_{0}:=I_{0}\left( \lambda ,h\right) =\int\limits_{s_{n}}^{s_{n-1}}\mathcal{C%
}_{n,\lambda }\left( s\right) ds=\frac{1-e^{-\lambda h}}{\lambda },
\end{equation*}%
\begin{equation*}
I_{1,n}:=I_{1,n}\left( \lambda ,h\right)
=\int\limits_{s_{n}}^{s_{n-1}}\left( s_{n-1}-s\right) \left[ s^{2}-s\left(
s_{n-1}-s\right) \right] \mathcal{C}_{n,\lambda }\left( s\right) ds,
\end{equation*}
\begin{equation*}
A_{1,n}:=A_{1,n}\left( \lambda ,h\right) =\frac{2}{I_{0}}\int%
\limits_{s_{n}}^{s_{n-1}}\left( s^{2}-2s\left( s_{n-1}-s\right) \right)
\mathcal{C}_{n,\lambda }\left( s\right) ds,
\end{equation*}%
\begin{equation*}
A_{2,n}:=A_{2,n}\left( \lambda ,h\right) =\frac{2}{I_{0}}\int%
\limits_{s_{n}}^{s_{n-1}}s\mathcal{C}_{n,\lambda }\left( s\right) ds.
\end{equation*}
In  equation (\ref{4.5}) the tail function $V$ is also
unknown. However, we observe that
\begin{equation}
\frac{\left\vert I_{1,n}\left( \lambda ,h\right) \right\vert }{I_{0}\left(
\lambda ,h\right) }\leq \frac{4\overline{s}^{2}}{\lambda }~\text{ for }%
\lambda h\geq 1.  \label{4.6}
\end{equation}%
Equation (\ref{4.6}) means that by taking $\lambda \gg 1,$ we mitigate
the influence of the nonlinear term with $\left( \nabla q_{n}\right)
^{2}$ in (\ref{4.5}). To solve  system (\ref{4.5})--(\ref{eq:q7}), we use
following algorithm: 

\subsubsection*{Globally convergent algorithm}

\begin{itemize}
\item Initialization: set $q_{0}\equiv 0$ and compute the first tail
  function $V_{0}$ as described in section 2.9 of \cite{BK} and
  \cite{BK2}.

\item For $n=1,\,2,\,\ldots ,\,N$

\begin{enumerate}
\item Set $q_{n,\,0}=q_{n-1}$, $V_{n,\,1}=V_{n-1}$

\item For $i=1,\,2,\,\ldots ,\,m_{n}$

\begin{itemize}
\item Find $q_{n,\,i}$ by solving (\ref{4.5})--(\ref{eq:q7}) with $%
V_{n}:=V_{n,\,i}$.

\item Compute $v_{n,i}=-hq_{n,\,i} - h \sum_{j=0}^{n-1} q_j + V_{n,\,i}$.

\item Compute $a_{n,i}$ via  discretization of (\ref{3.10}) with $a:= a_{n,i}$
 and $v:=v_{n,i}$. Then
solve the forward problem (\ref{2.1})--(\ref{2.2}) with the new
computed coefficient $a:= a_{n,i}$, compute $w:=w_{n,\,i}$ and update the tail $V_{n,\,i+1}$ by (\ref{2.A}).
\end{itemize}

\item Set $q_{n}=q_{n,\,m_{n}}$, $a_n = a_{n,m_{n}}$, $V_{n}=V_{n,\,m_{n+1}}$ and go to the next frequency interval $\left[ s_{n+1},\,s_{n}\right] $ if $n<N.$ If $n=N$, then stop.
\end{enumerate}
\end{itemize}

The stopping criteria for iterations $m_n$ and $n$ and step 3 in the above algorithm is derived computationally in  \cite{BTKF,NBKF}. The global
convergence theorem was proven in \cite{BK, BK2}. 

\section{Finite element method for  reconstruction}
\label{sec:fem}

In this section we explain how we can reconstruct the function $a(x)$
of the equation (\ref{2.1}) using the variational formulation of
equation (\ref{2.7}).  Suppose that the pair of functions $\left(
V_{n,i},q_{n,i}\right) $ at step 2 of the globally convergent
algorithm is computed. Then using  the Finite Difference discretization of (\ref{3.7}) we can compute the function $v_{n,i}\left(x\right)$ as
\begin{equation}
v_{n,i}\left( x\right) =-hq_{n,i}\left( x\right)
-h\sum\limits_{j=0}^{n-1}q_{j}\left( x\right) +V_{n,i}\left( x\right)~\text{
}x\in \Omega.  \label{3.104}
\end{equation}
Using (\ref{transform})  we can get
\begin{equation}
v_{n,i}\left( x\right) =\frac{\ln w_{a_{n,i}}\left( x,s_{n}\right) }{%
s_{n}^{2}},  \label{3.105}
\end{equation}
and thus
\begin{equation*}
w_{a_{n,i}}\left(x\right) = \mathrm{e}^{s_{n}^{2}v_{n,i}\left( x\right)}.
\end{equation*}
Here, the function $w_{a_{n,i}}\left( x,s_{n}\right) $ is the solution of
the following analog of the problem (\ref{2.7}), (\ref{2.8})%
\begin{equation}
\Delta w_{a_{n,i}}-s_{n}^{2}a_{n,i}\left( x\right) w_{a_{n,i}}=0\text{ in }%
\Omega ,  \label{3.106}
\end{equation}%
\begin{equation}
\partial _{n}w_{a_{n,i}}\mid _{\partial \Omega }=f_{n,i}\left( x\right) ,
\label{3.107}
\end{equation}%
where
\begin{equation*}
f_{n,i}\left( x\right) =\partial _{n}\mathrm{e}^{ s_{n}^{2}v_{n,i}\left(
x\right)}\text{ for }x\in \partial \Omega.
\end{equation*}%

To find $a_{n,i}$ from (\ref{3.106}), we will use the finite
element method for the problem (\ref{3.106})--(\ref{3.107}).  We
introduce the finite element trial space $V_h$, defined by
\begin{equation}
V_h := \{ u \in  H^1(\Omega): u|_{K} \in  P_1(K),  \forall K \in K_h \}, \nonumber
\end{equation}
where $P_1(K)$ denotes the set of linear functions on the element $K$
of the finite element mesh $K_h$.
Hence, the finite element space $V_h$
consists of continuous piecewise linear functions in space.
To approximate functions $a_{n,i}$  we introduce space of piecewise-linear functions  $C_h$ defined by
\begin{equation}
C_h := \{ u \in  H^1(\Omega): u|_{K} \in  P_1(K),  \forall K \in K_h \}, \nonumber
\end{equation}
%where $P_0(K)$ is the piecewise constant function defined in the vertices of th%e element $K$ of the mesh $K_h$.
Let us define a $L_2$ inner product
\begin{equation*}
(\alpha ,\beta )=\int_\Omega\alpha\beta \, dx. 
\end{equation*}
Then the finite element formulation for (\ref{3.106})-(\ref{3.107})  reads:
Find $a_{n,i} \in C_h,  w_{a_{n,i}} \in V_h$ such that   for all $v \in V_h$
\begin{equation}
( a_{n,i}  w_{a_{n,i}}, v) =  -\frac{1}{s_{n}^{2}}(\nabla w_{a_{n,i}}, \nabla  v ) + \frac{1}{s_{n}^{2}} ( f_{n,i}, v )_{\partial \Omega}.  \label{3.107_1}
\end{equation}

We expand $w_{a_{n,i}}$ in terms of the standard continuous piecewise
linear functions $\{\varphi_k\}_{k=1}^N$ in space as 
\begin{equation}\label{new}
w_{a_{n,i}}=\sum_{k=1}^N {w_{a_{n,i}}}_k \varphi_k(x),
\end{equation}
 where ${w_{a_{n,i}}}_k $ are the discrete nodal values of the
already computed functions  $v_{n,i}$ at step 2 of the globally convergent algorithm  with the nodal values ${v_{n,i}}_k$ such that
\begin{equation*}
 {w_{a_{n,i}}}_k  = \mathrm{e}^{ s_{n}^{2}{v_{n,i}}_k\left(x\right)}~ \forall x \in \Omega.
\end{equation*}
We substitute 
expansion (\ref{new}) in the variational formulation (\ref{3.107_1}) with $v(x) =
\varphi_j(x)$, and obtain the following system of discrete equations
\begin{equation}
\sum_{k,j=1}^N {a_{n,i}}_k ( {w_{a_{n,i}}}_k~ \varphi_k, \varphi_j) =  -\frac{1}{s_{n}^{2}} \sum_{k,j=1}^N  {w_{a_{n,i}}}_k( \nabla \varphi_k, \nabla  \varphi_j ) +
\frac{1}{s_{n}^{2}} \sum_{j=1}^N  ( f_{n,i}, \varphi_j)_{\partial \Omega}.  \label{3.108}
\end{equation}
The system (\ref{3.108}) can be rewritten in the matrix form for the unknown $a_{n,i}$ and known $ w_{a_{n,i}}$ as
\begin{equation}
M  a_{n,i} =  -\frac{1}{s_{n}^{2}}  G  w_{a_{n,i}}  + \frac{1}{s_{n}^{2}} F.
 \label{3.108_1}
\end{equation}
Here, $M$ is the block mass matrix in space, $G$ is the stiffness
matrix corresponding to the gradient term, $F$ is the load vector.  At
the element level the matrix entries in (\ref{3.108_1}) are explicitly
given by:
 \begin{eqnarray}
  M_{k,j}^{K} & = &   ( {w_{a_{n,i}}}_k ~\varphi_k, \varphi_j)_K, \\
  G_{k,j}^{K} & =  & (\nabla \varphi_k, \nabla  \varphi_j)_K,\\
  F_{j}^{K}&=& (f_{n,i}, \varphi_j)_{K}.
 \end{eqnarray}

To obtain an explicit scheme for the computation of the coefficients
$a_{n,i}$, we approximate $M$ by the lumped mass matrix $M^{L}$ in
space. This matrix is obtained as the diagonal approximation of the
mass matrix $M$: diagonal elements of $M^L$ are obtained as the row
sum of elements in $M$. Thus, we get the following equation for the
explicit computation of the function $a_{n,i}$ in (\ref{3.106}):
\begin{equation}
a_{n,i} = -\frac{1}{s_{n}^{2}}  (M^{L})^{-1}  G w_{a_{n,i}}  + \frac{1}{s_{n}^{2}}  (M^{L})^{-1}  F.
  \label{3.109}
\end{equation}

\begin{figure}[tbp]
\begin{center}
\begin{tabular}{ccc}
{\includegraphics[scale=0.25,clip=]{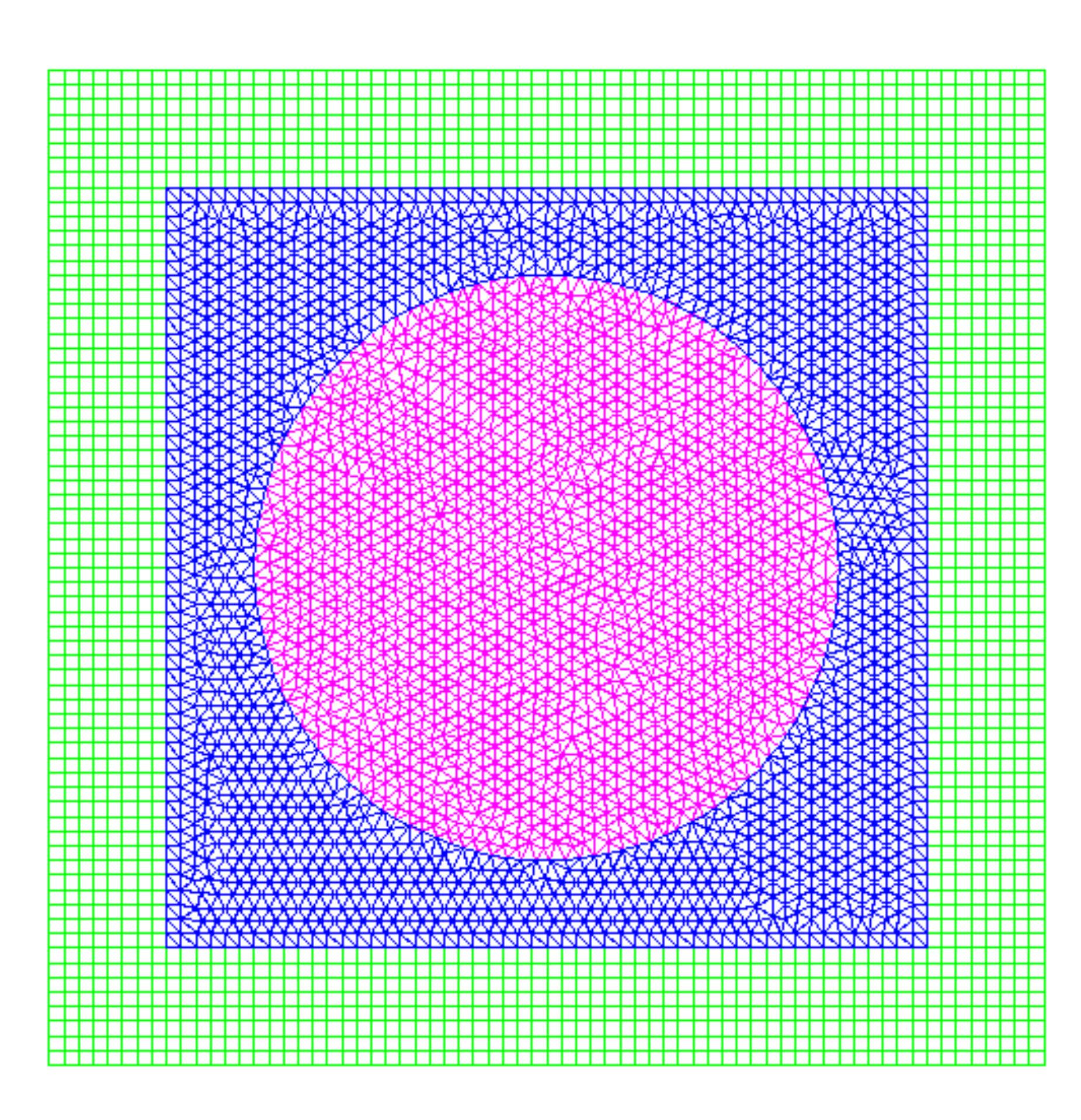}} & 
{\includegraphics[scale=0.25,clip=]{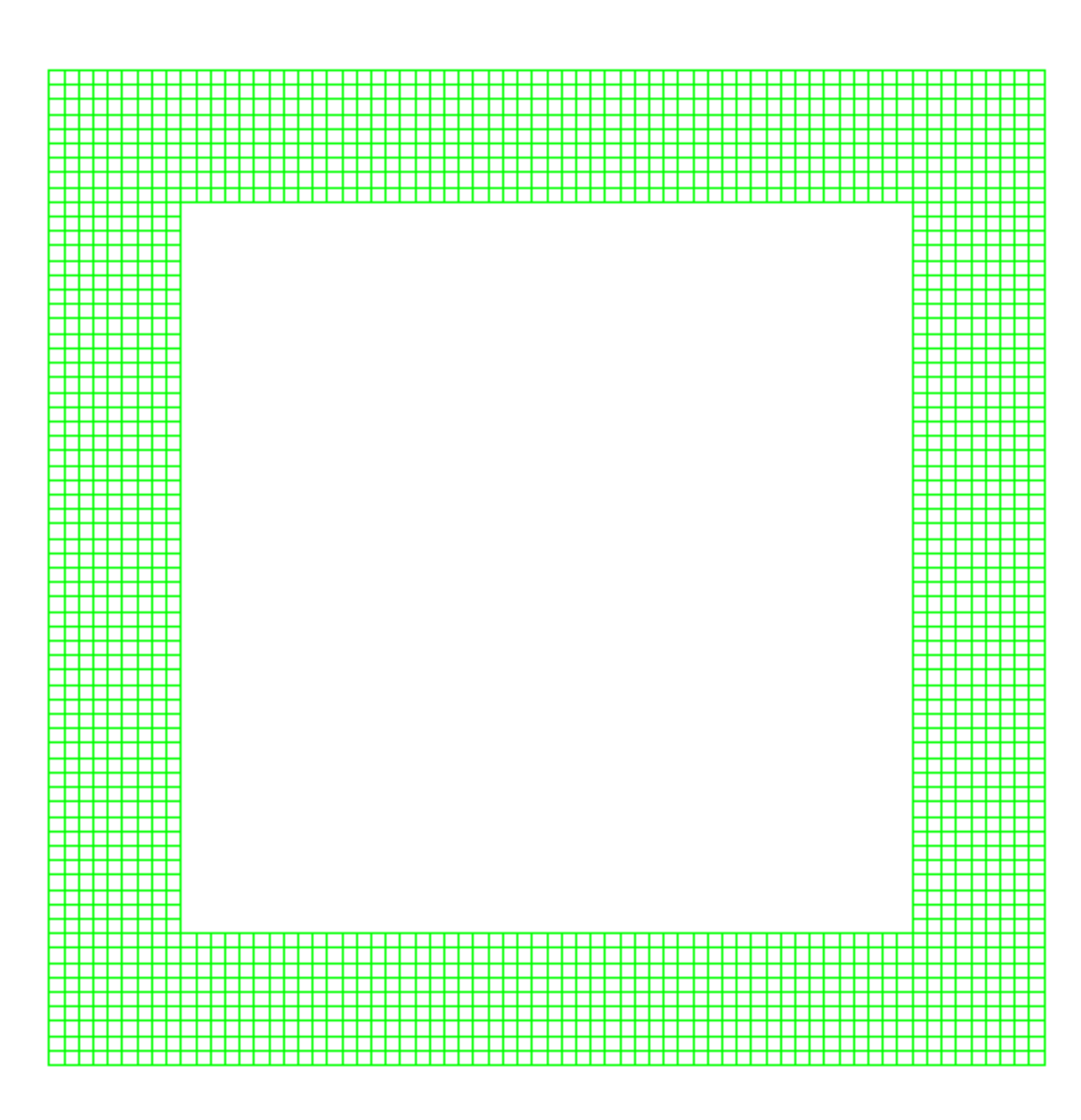}} & 
{\includegraphics[scale=0.25,clip=]{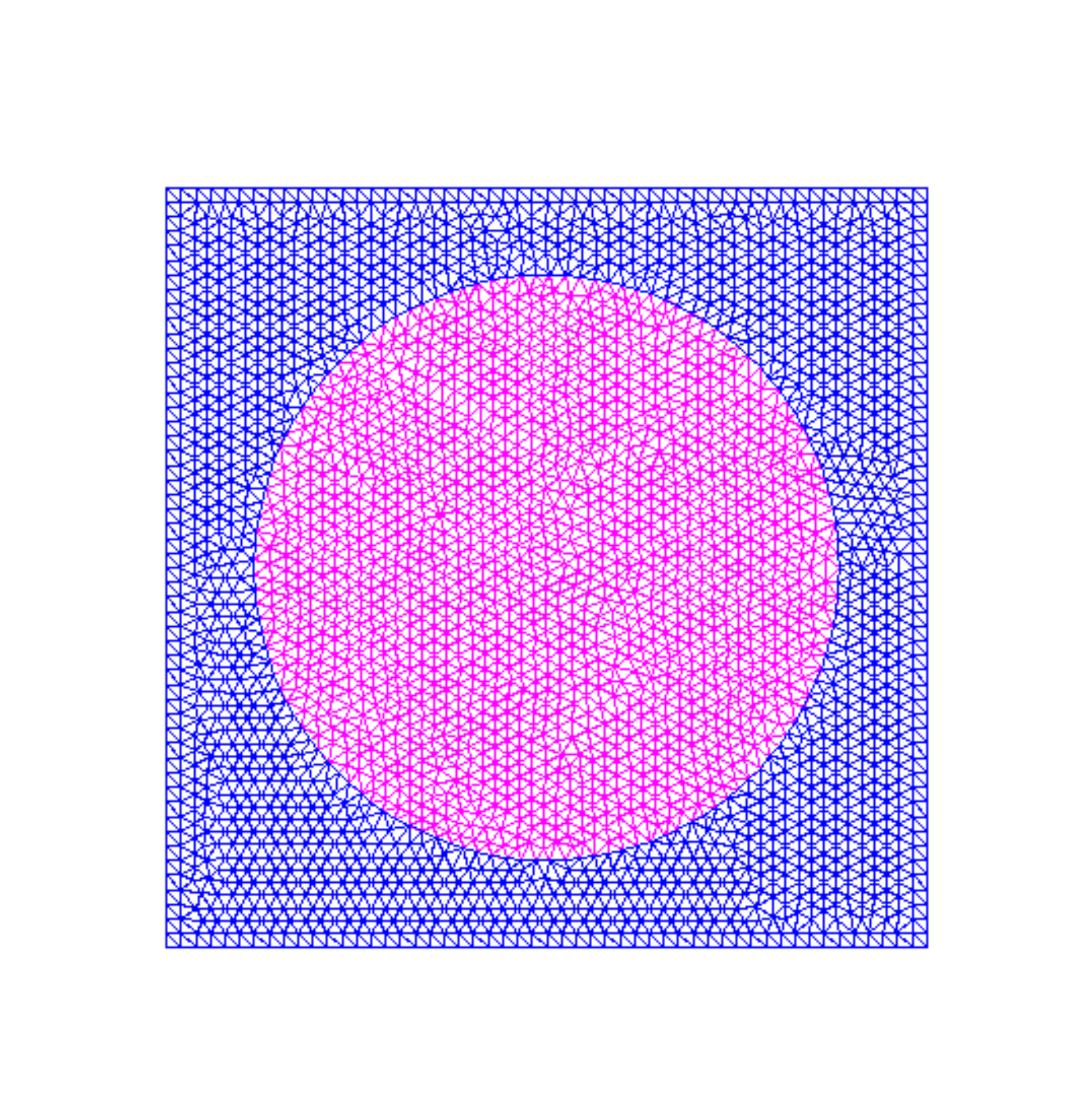}} \\ 
a) $ G = G_{FEM} \cup G_{FDM}$ & b) $ G_{FDM}$ & 
c) $G_{FEM}= \Omega$
\end{tabular}
\end{center}
\caption{ \emph{\ a) Geometry of the hybrid mesh. This is a
    combination of the quadrilateral finite difference mesh in the subdomain
    $G_{FDM}$ presented on b),  and the finite element mesh
    in the inner domain $G_{FEM}= \Omega$ presented on c). The solution of the inverse problem is computed
    in $G_{FEM}= \Omega$. We use software package WavES \cite{waves} to compute hybrid solution on these meshes. }}
\label{fig:F1}
\end{figure}

\begin{figure}[tbp]
\begin{center}
\begin{tabular}{ccc}
{\includegraphics[scale=0.2,clip=]{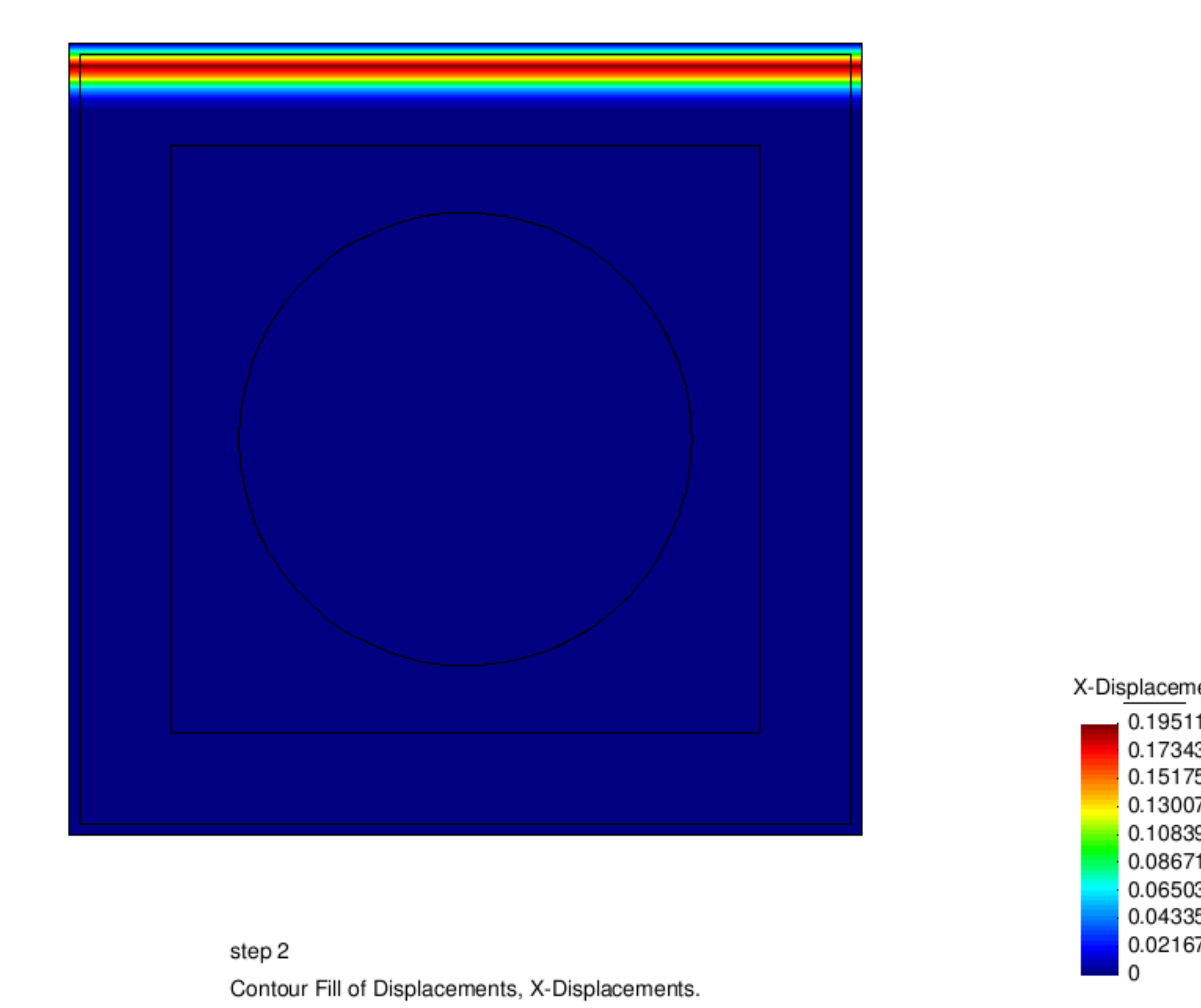}} & {
\includegraphics[scale=0.2,clip=]{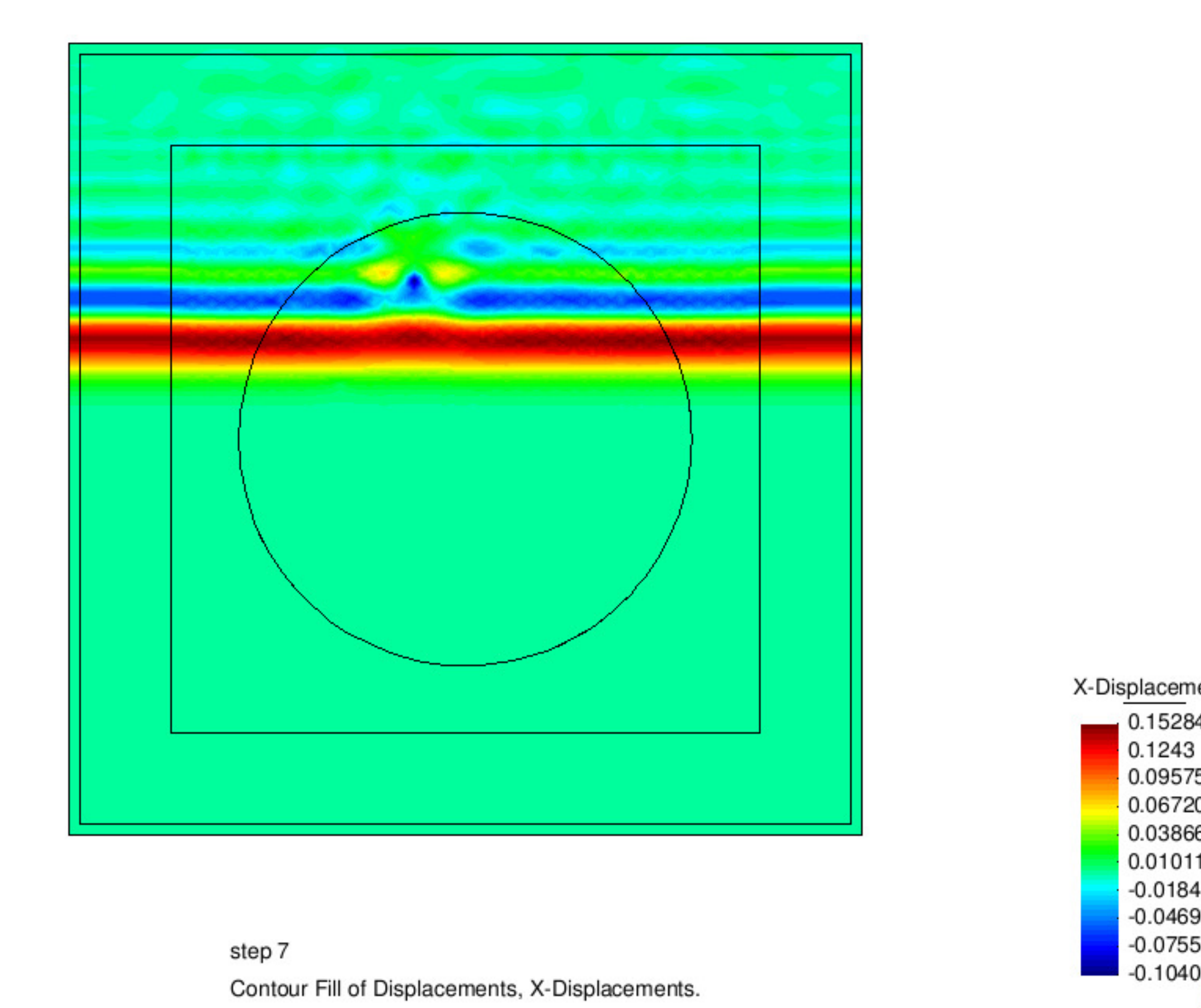}} & {
\includegraphics[scale=0.2,clip=]{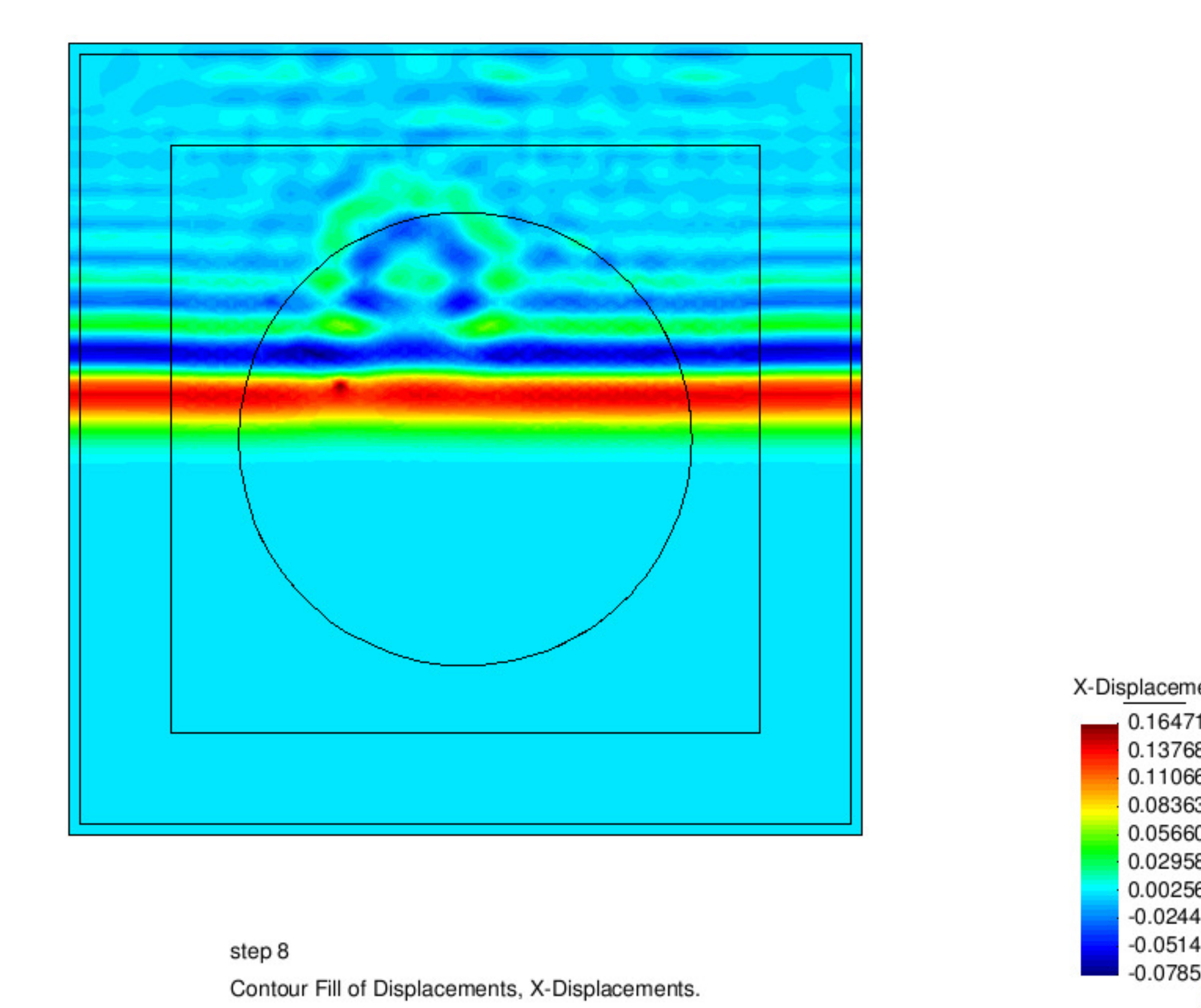}} \\ 
a) $t=0.2$ & b) $t=0.7$ & c) $t=0.8$ \\ 
{\includegraphics[scale=0.2,clip=]{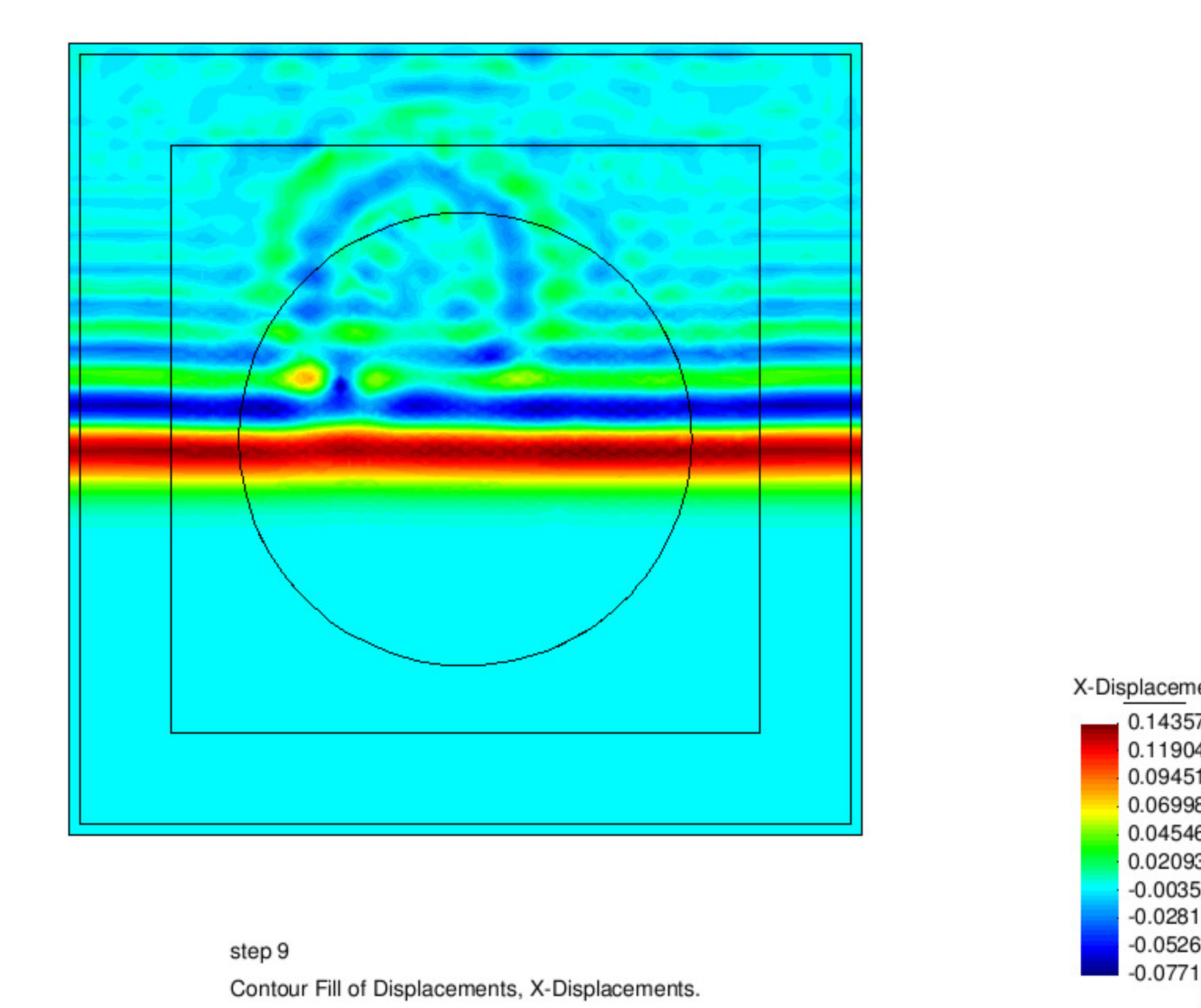}} & {
\includegraphics[scale=0.2,clip=]{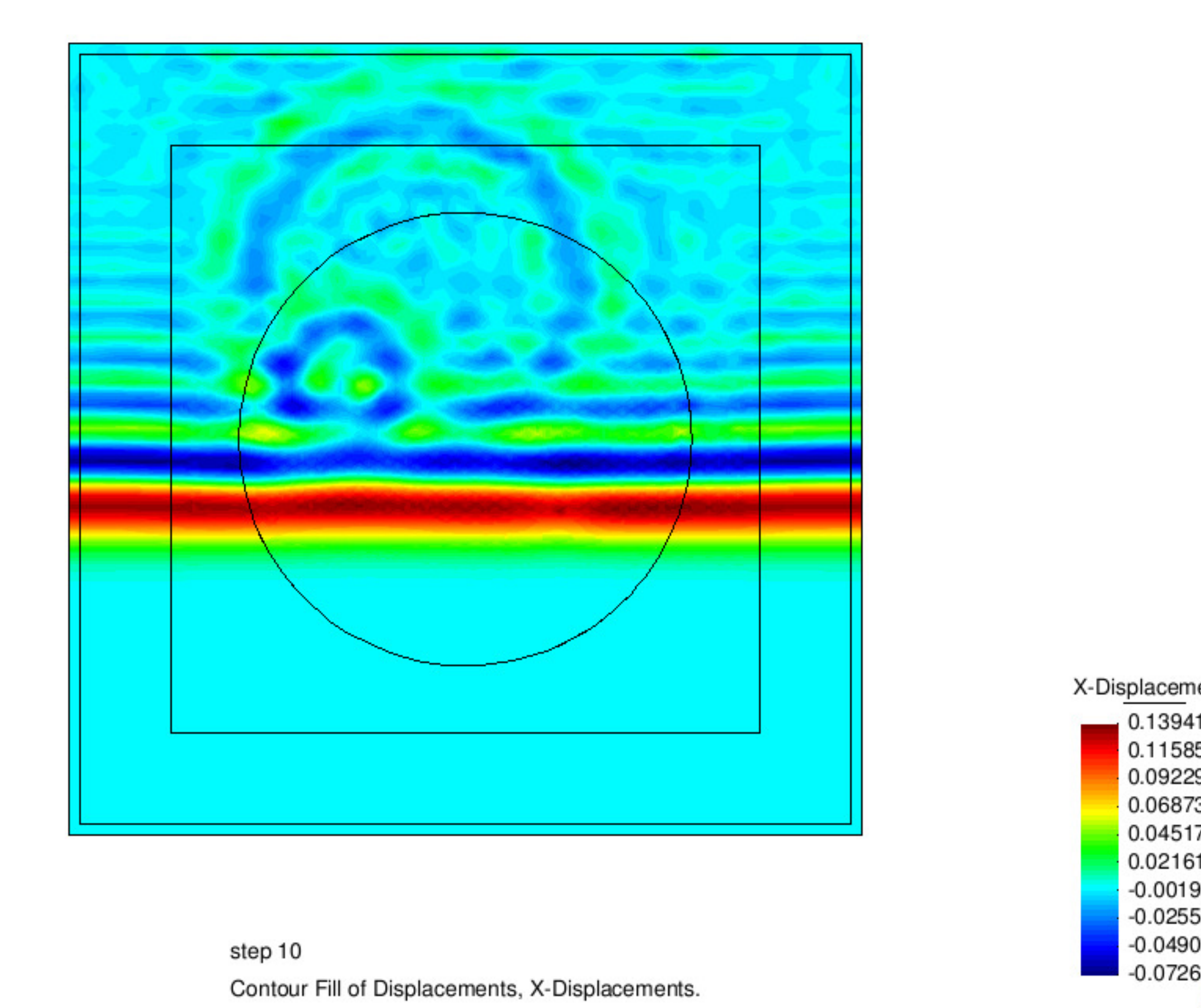}} & {
\includegraphics[scale=0.2,clip=]{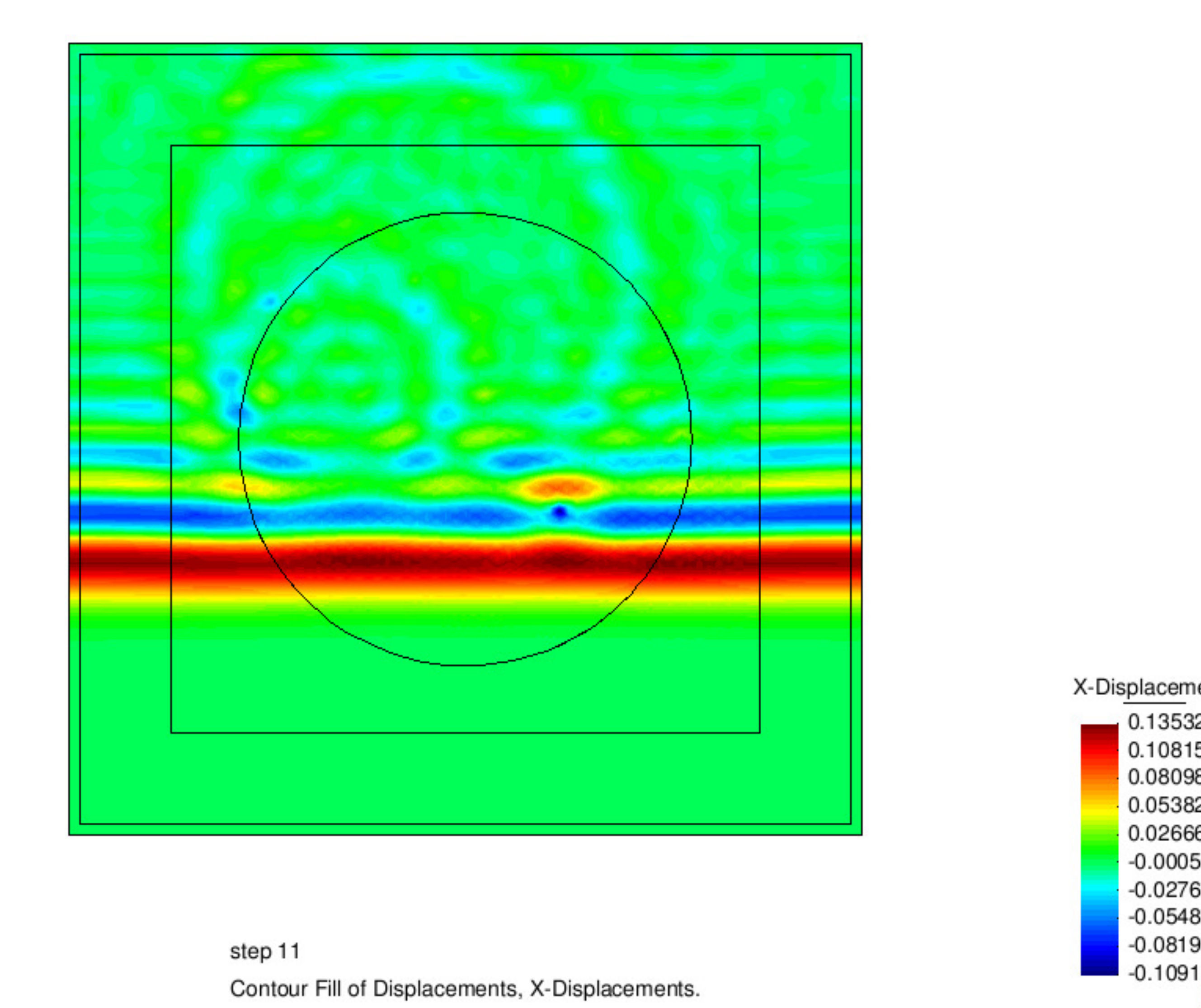}} \\ 
d) $t=0.9$ & e) $t=1.0$ & f) $t=1.1$ \\ 
{\includegraphics[scale=0.2,clip=]{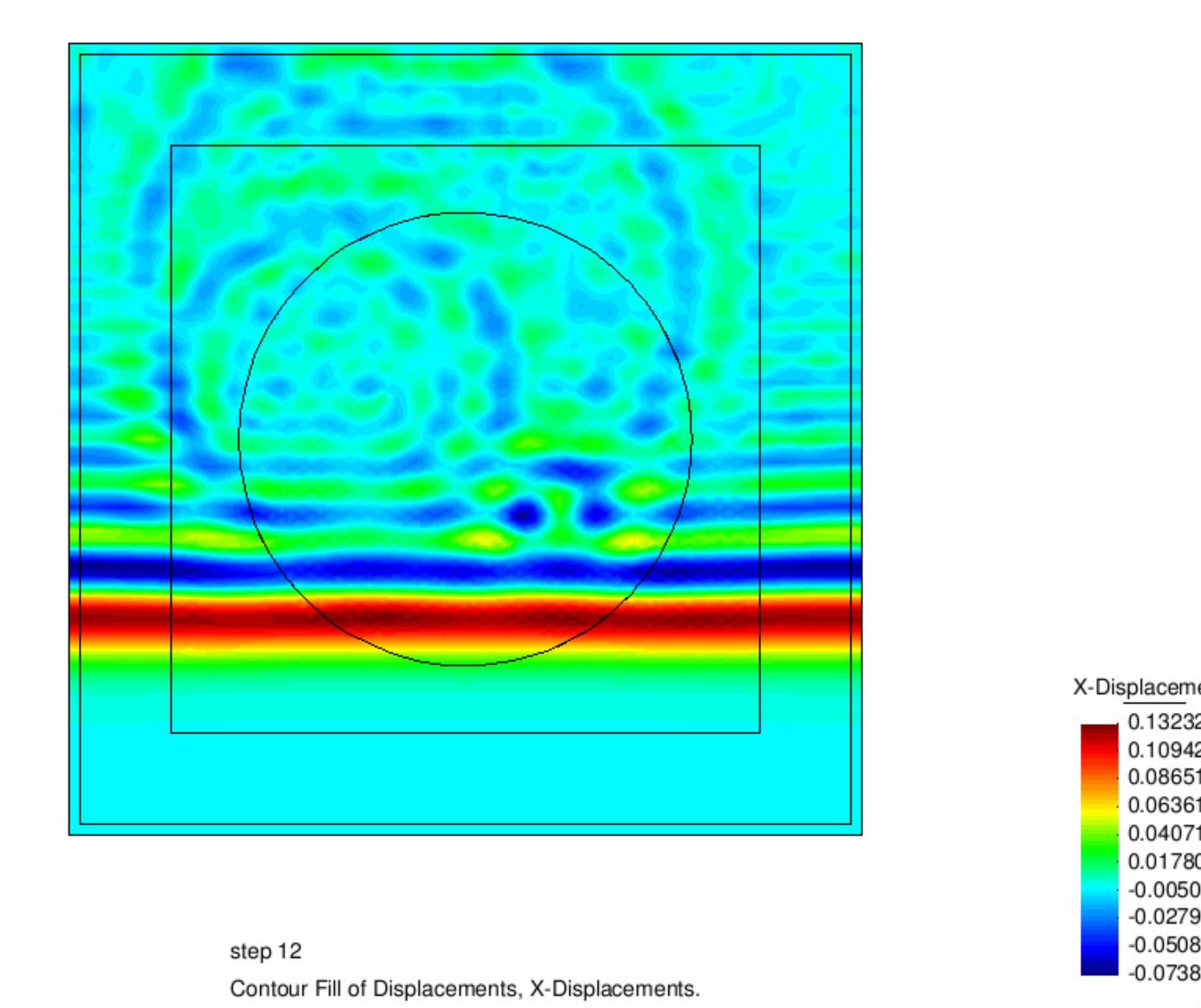}} & {
\includegraphics[scale=0.2,clip=]{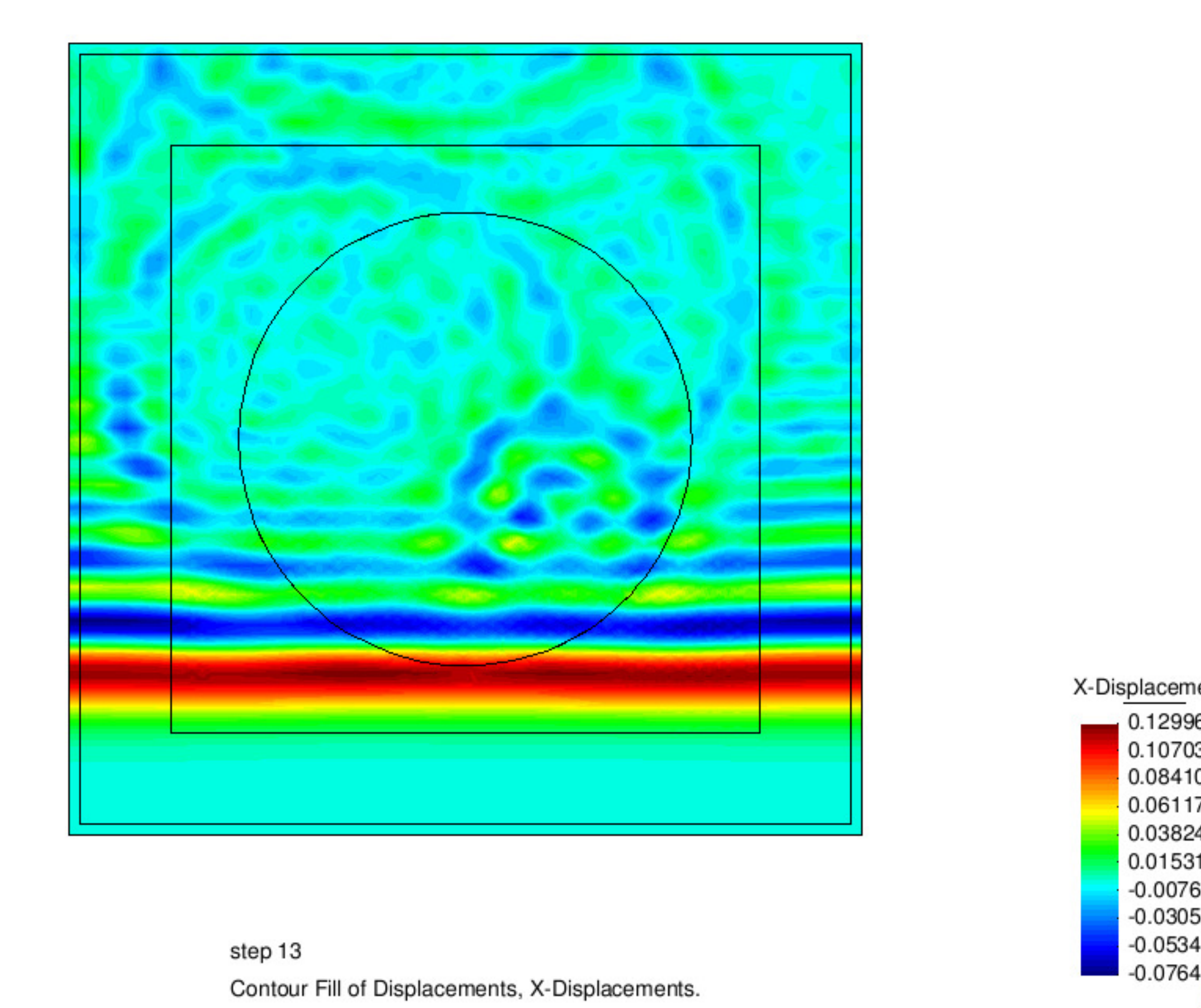}} & {
\includegraphics[scale=0.2,clip=]{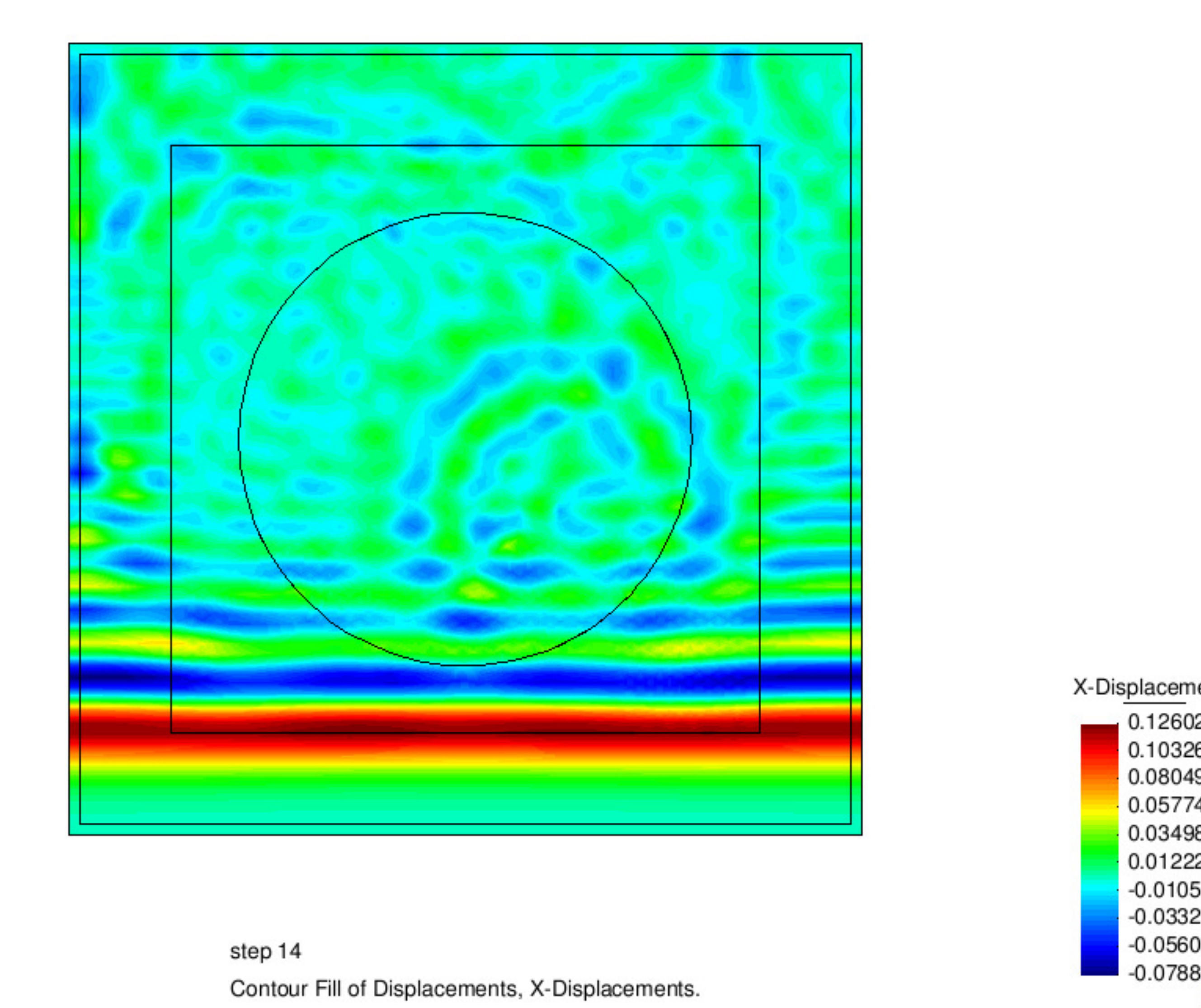}} \\ 
g) $t=1.2$ & h) $t=1.3$ & i) $t=1.4$ \\ 
{\includegraphics[scale=0.2,clip=]{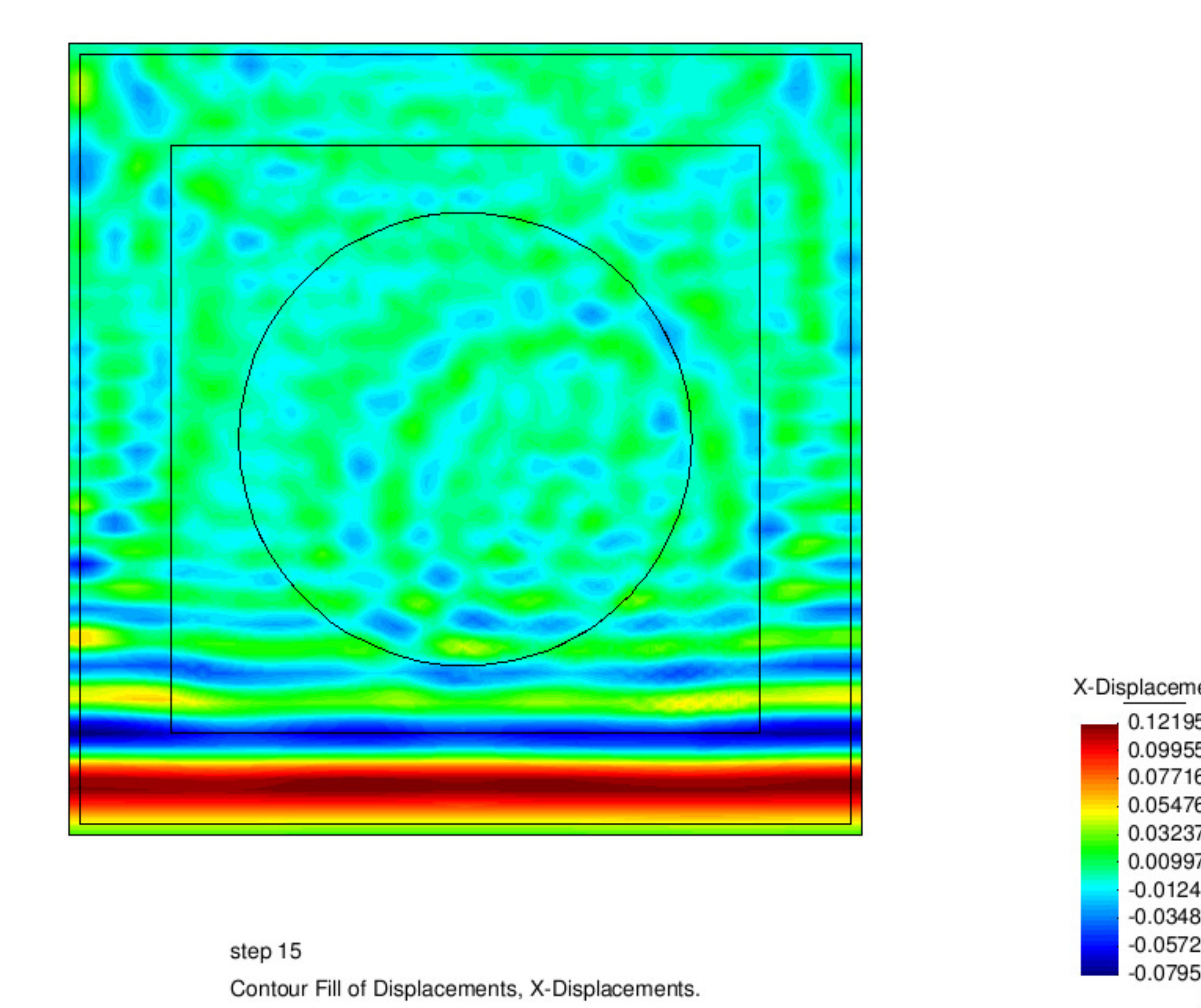}} & {
\includegraphics[scale=0.2,clip=]{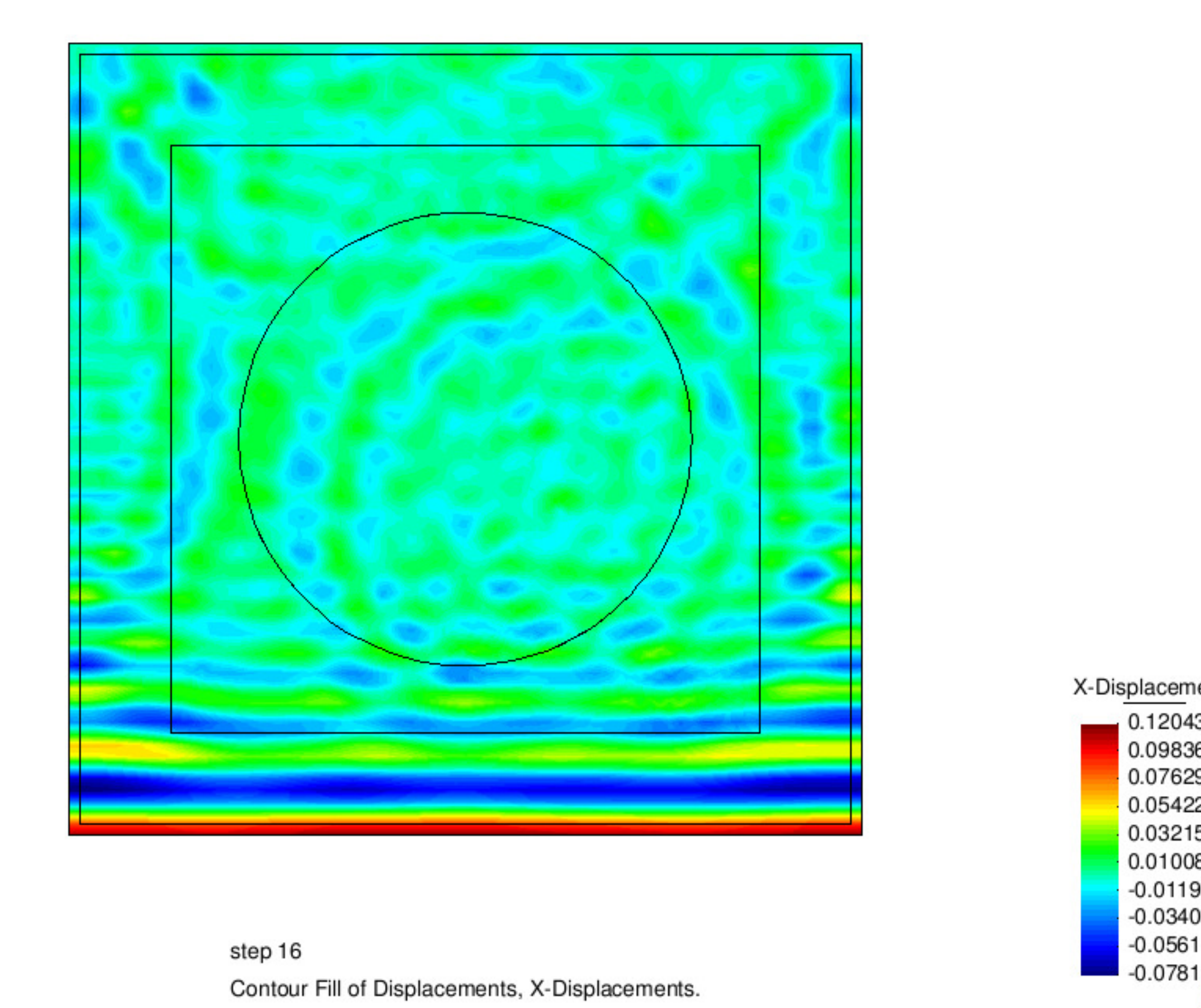}} & {
\includegraphics[scale=0.2,clip=]{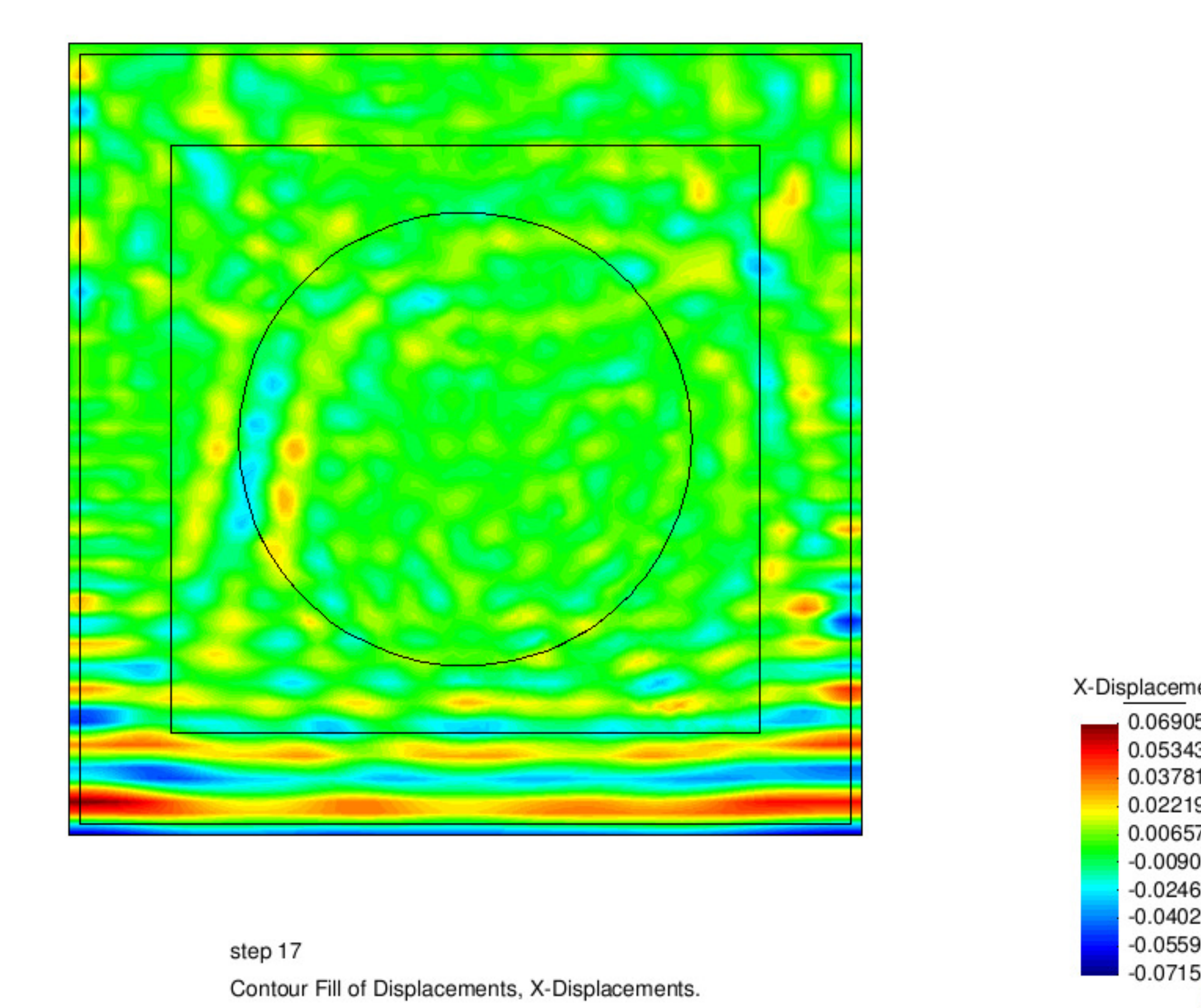}} \\ 
j) $t=1.5$ & k) $t=1.6$ & l) $t=1.7$ \\ 
&  & 
\end{tabular}%
\end{center}
\caption{{\protect\small \emph{Isosurfaces of the computed solution
      $u(x,t)$ of the wave equation in $G$ at different times $t$ with
      the plane wave initialized at the front boundary of the domain
      $G$.   Test was computed
      in time $t=[0,2]$ with time step $\protect\tau =
      0.001$. Software package WavES \cite{waves} is used for the numerical
      simulation of this solution. }}}
\label{fig:F3D_1}
\end{figure}

\begin{figure}[tbp]
\begin{center}
\begin{tabular}{ccc}
{\includegraphics[scale=0.15,clip=]{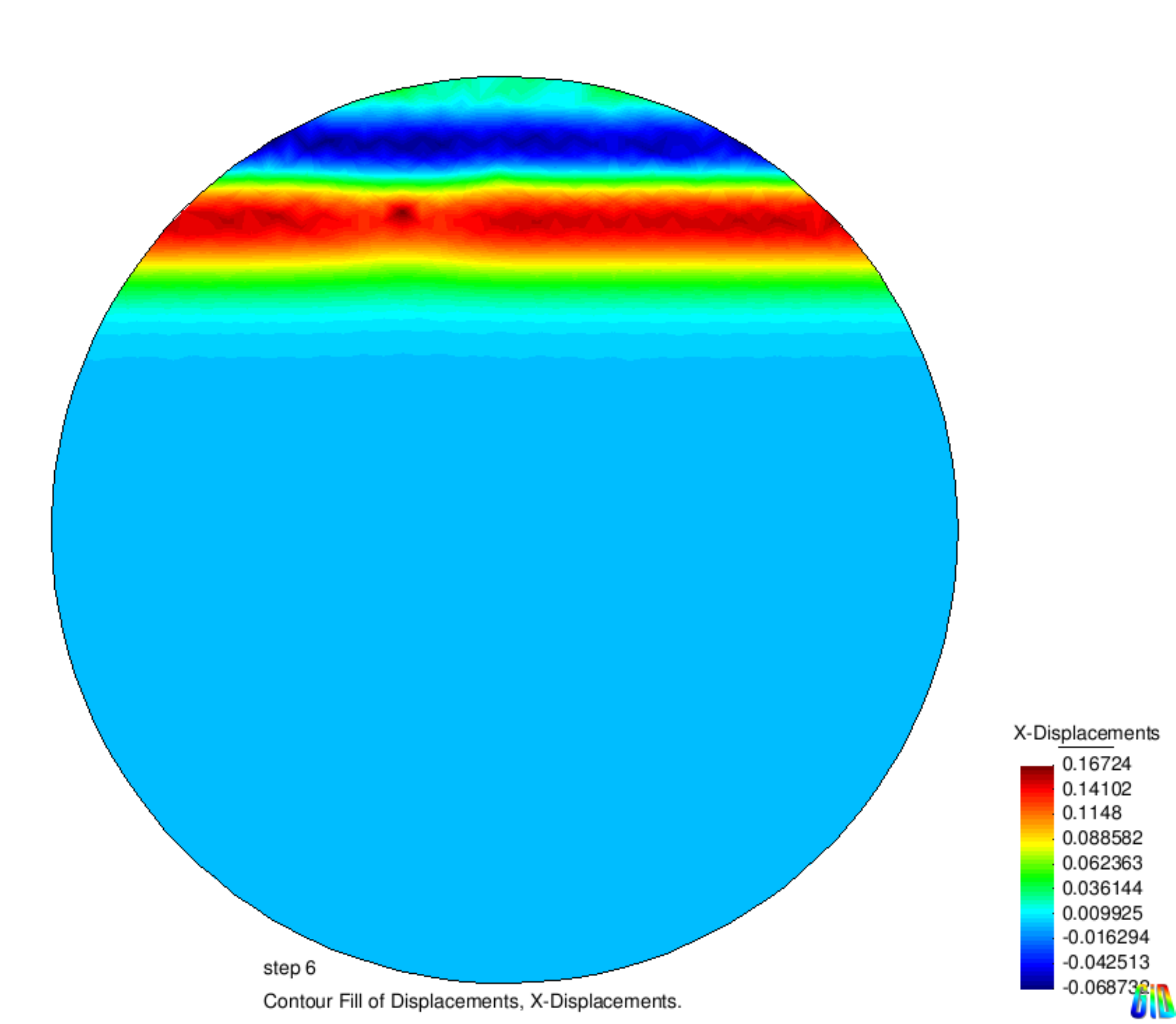}} & {%
\includegraphics[scale=0.15,clip=]{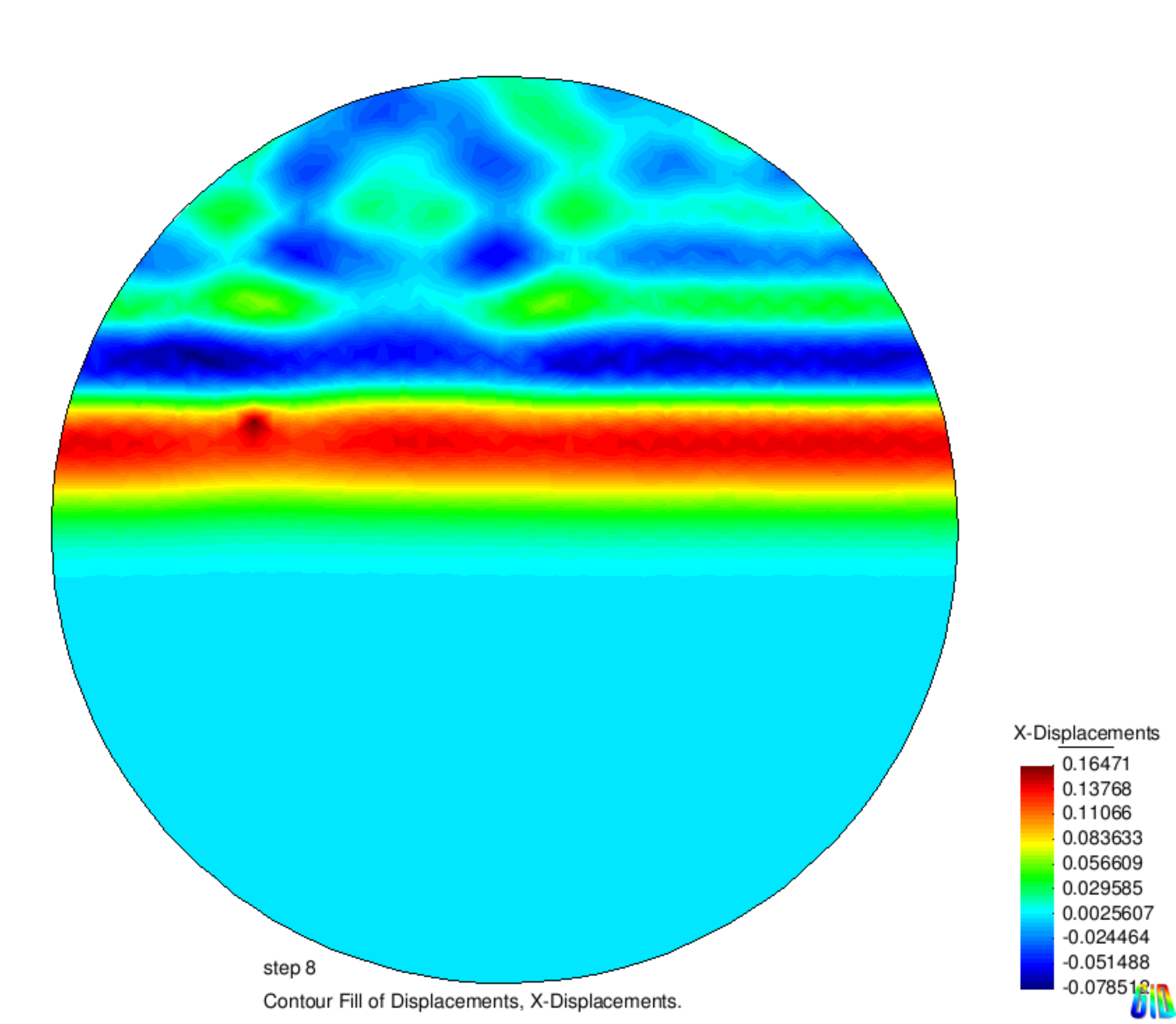}} & {%
\includegraphics[scale=0.15,clip=]{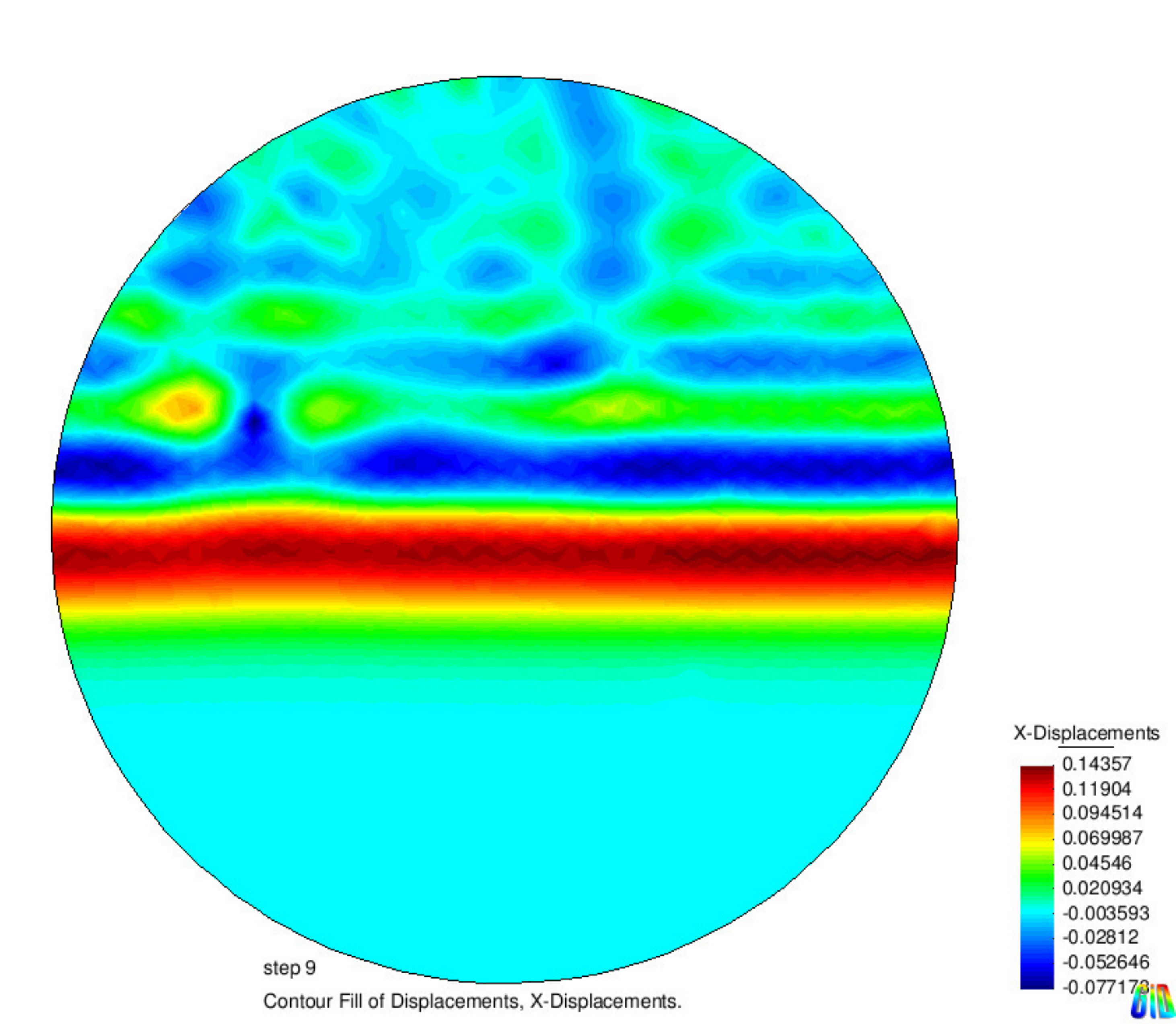}} \\ 
a) $t=0.6$ & b) $t=0.8$ & c) $t=0.9$ \\ 
{\includegraphics[scale=0.15,clip=]{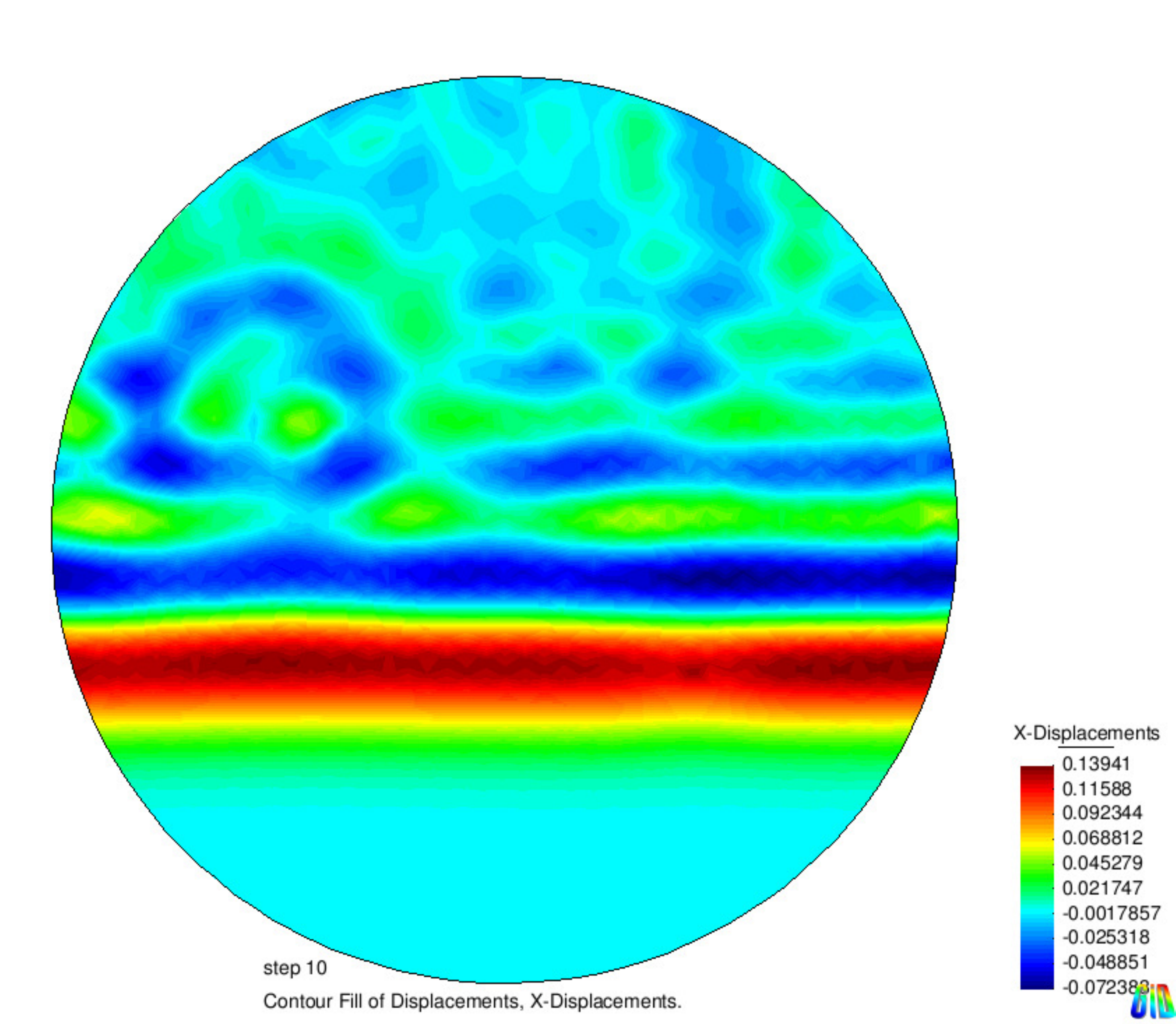}} & {%
\includegraphics[scale=0.15,clip=]{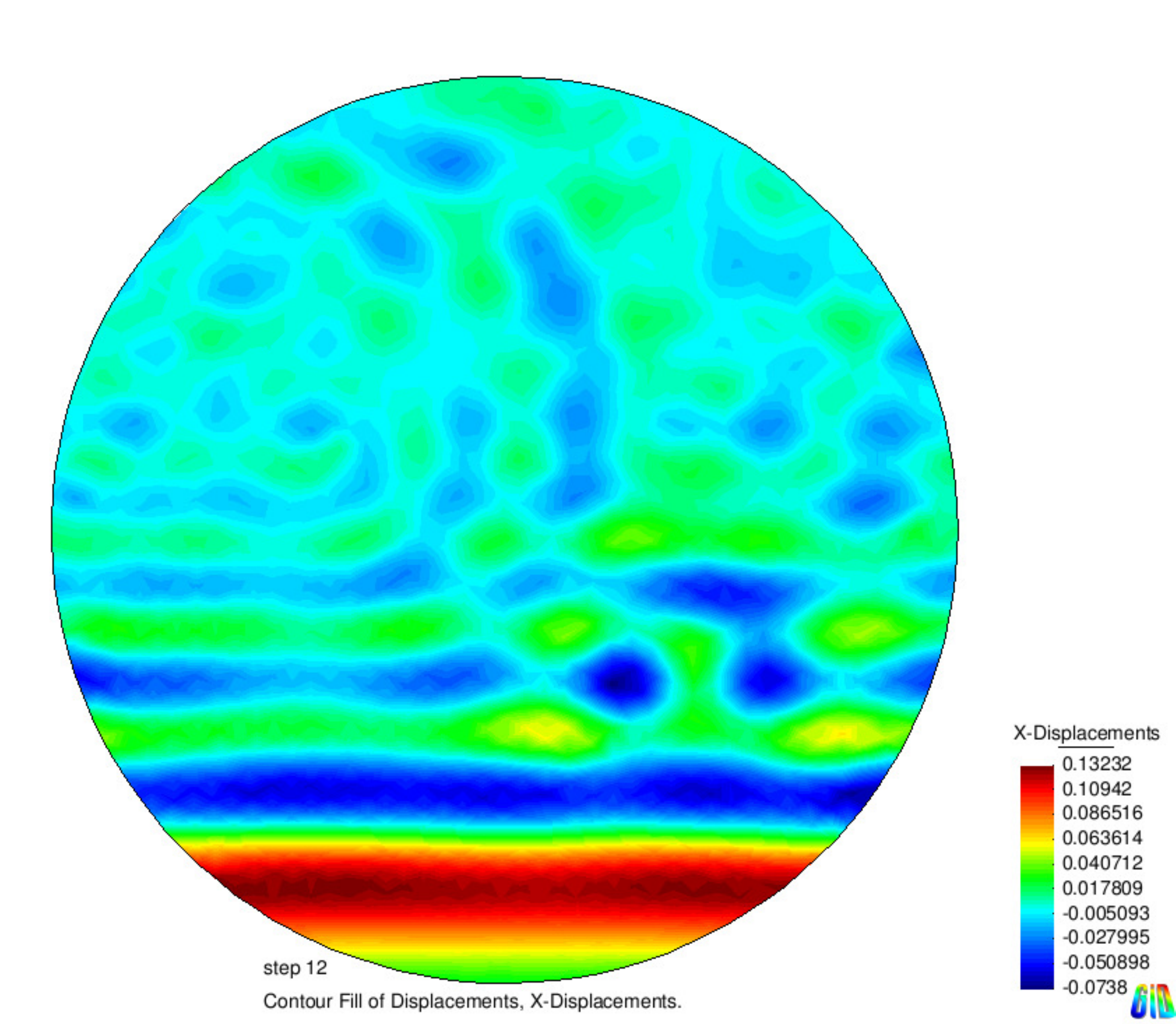}} & {%
\includegraphics[scale=0.15,clip=]{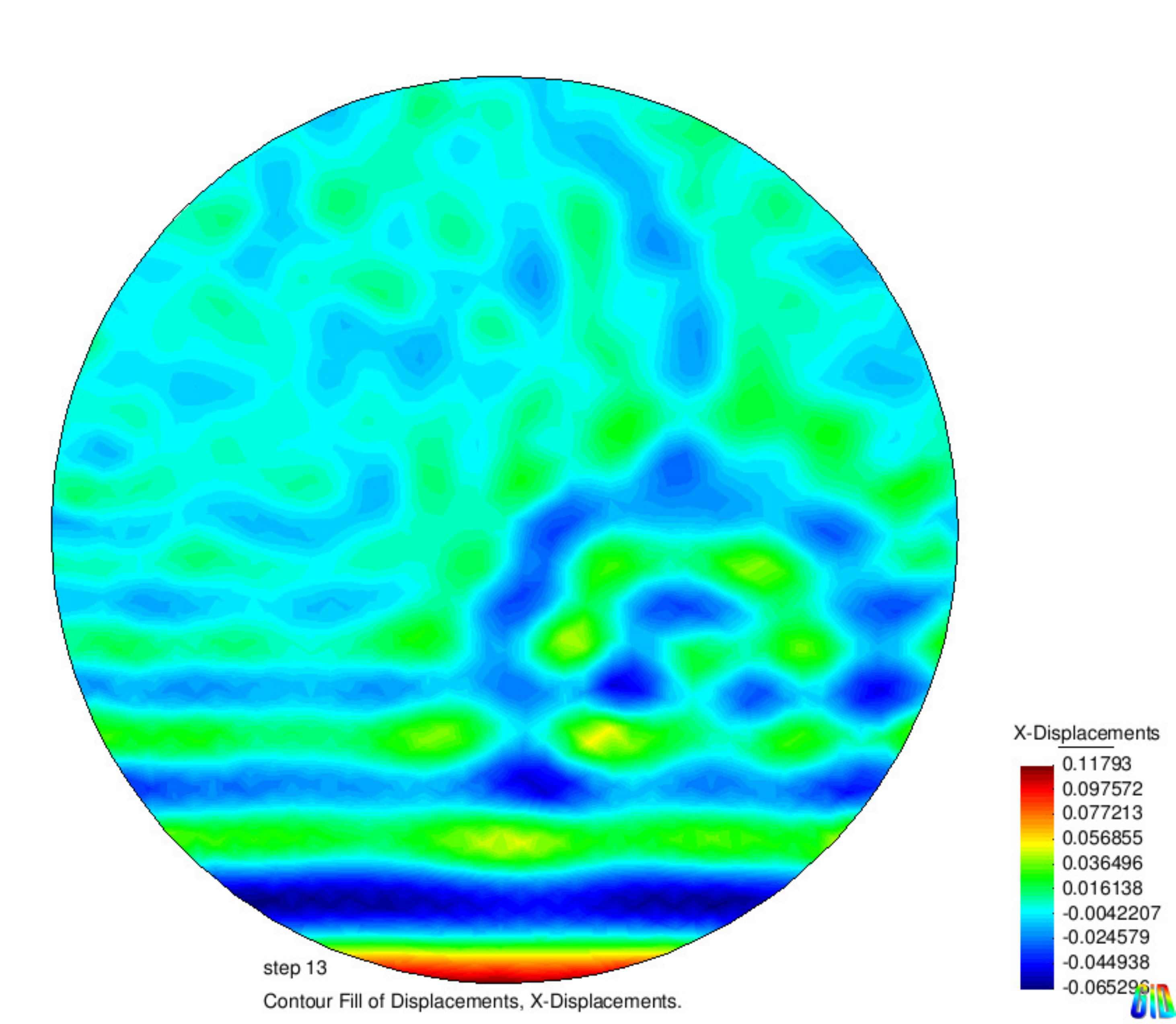}} \\ 
d) $t=1.0$ & e) $t=1.2$ & f) $t=1.3$ \\ 
{\includegraphics[scale=0.15,clip=]{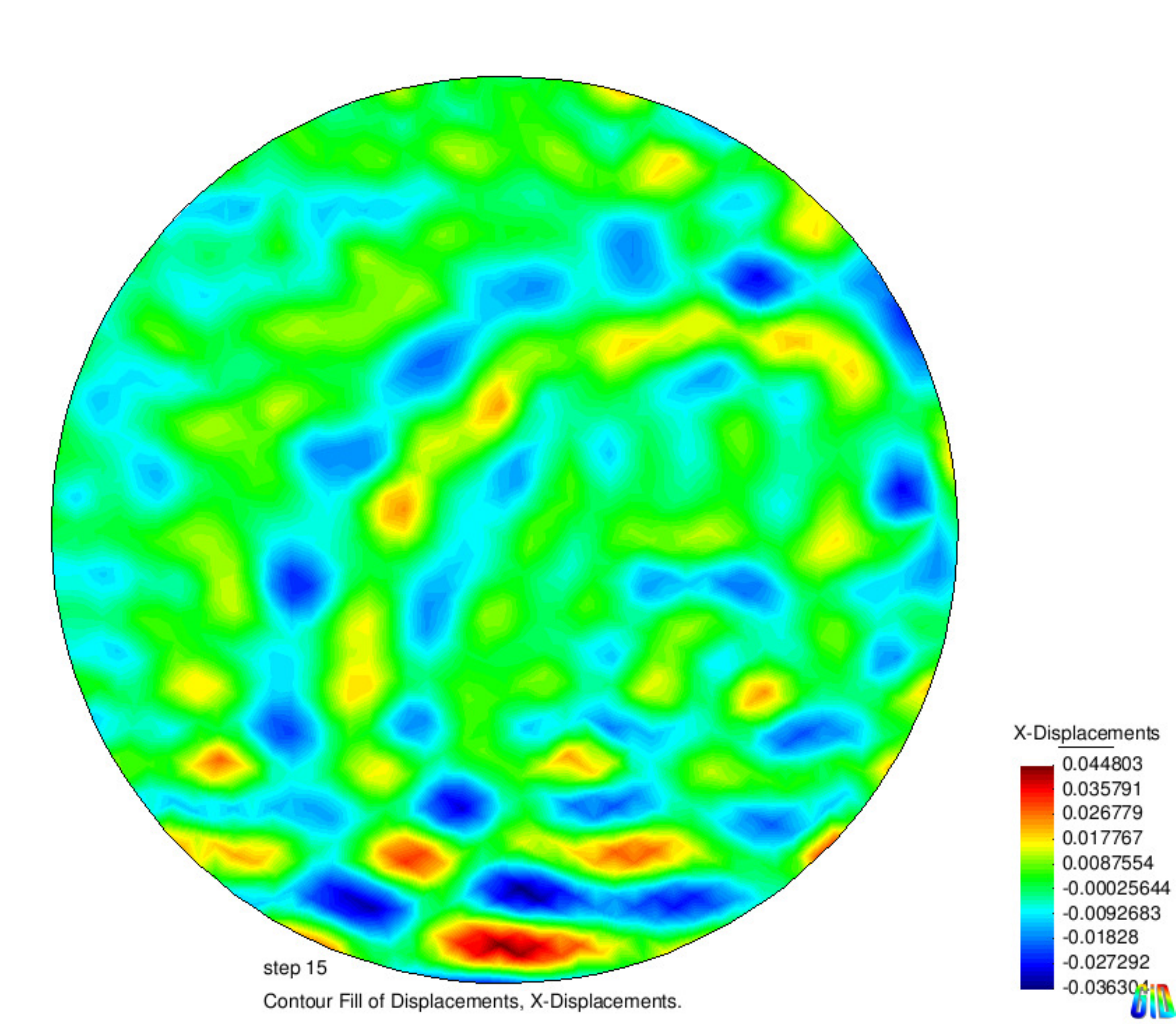}} & {%
\includegraphics[scale=0.15,clip=]{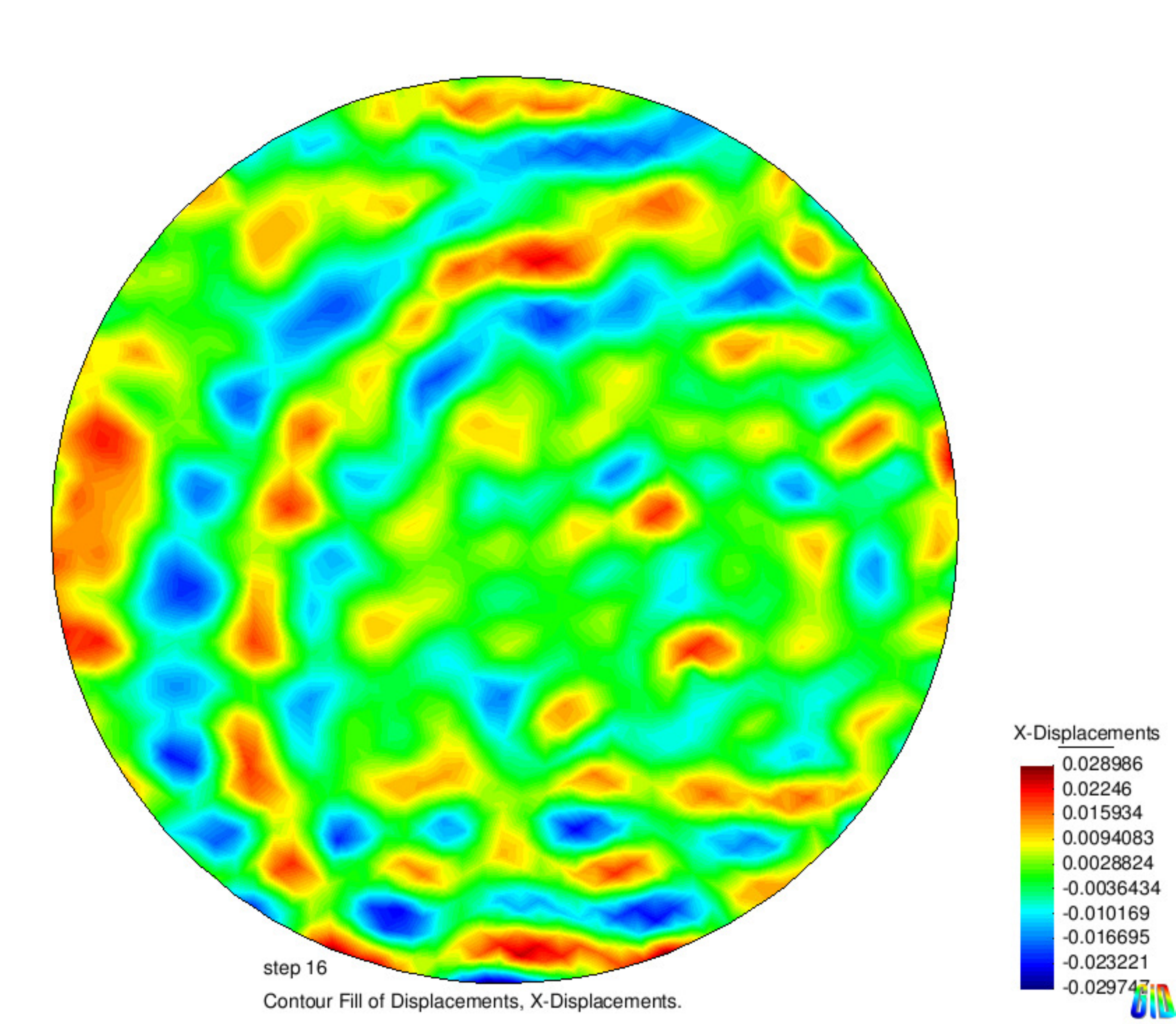}} & {%
\includegraphics[scale=0.15,clip=]{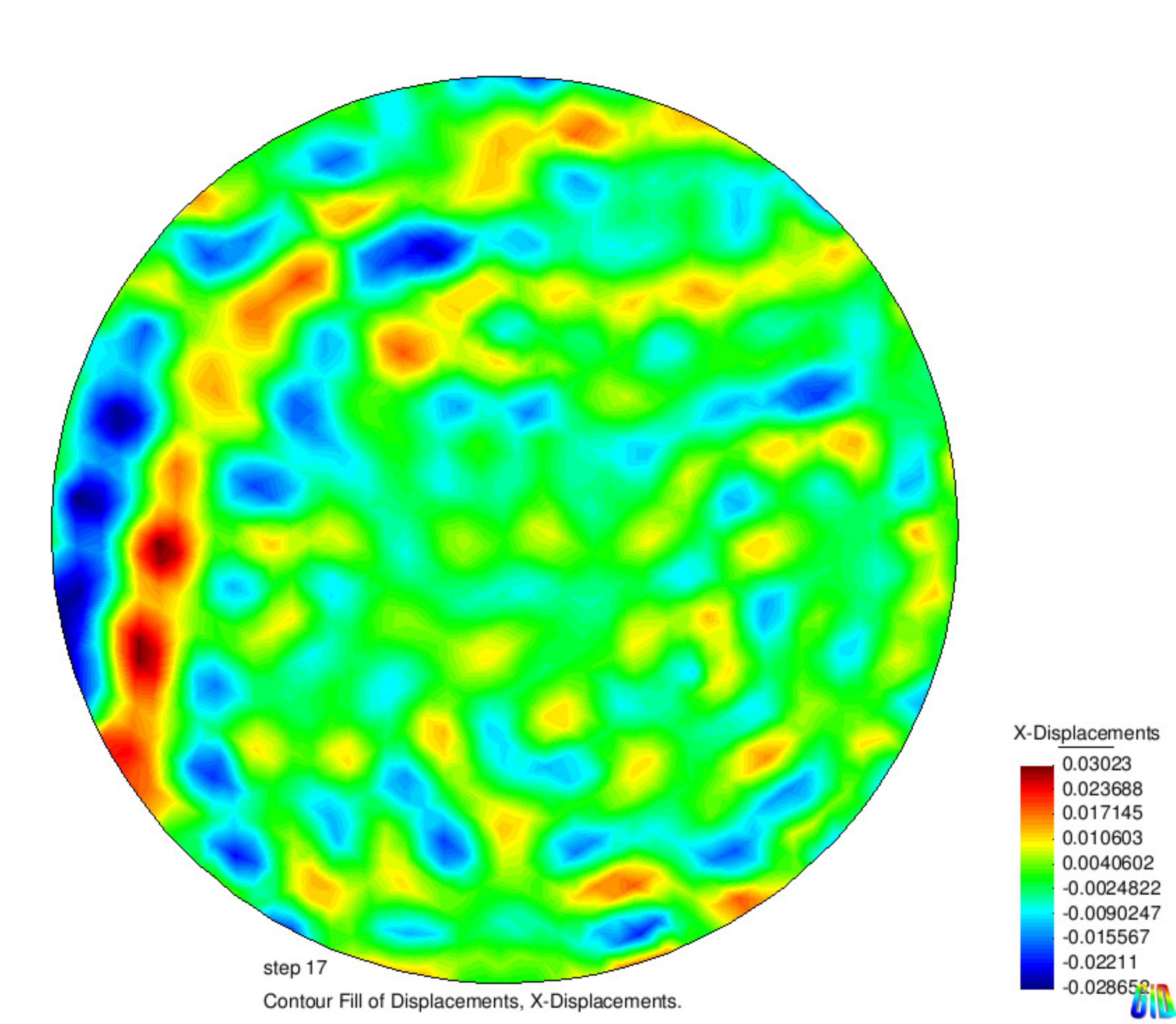}} \\ 
g) $t=1.5$ & h) $t=1.6$ & i) $t=1.7$ \\ 
&  & 
\end{tabular}%
\end{center}
\caption{{\protect\small \emph{Extracted isosurfaces of the computed
      solution $u(x,t)$ of Figure  \ref{fig:F3D_1} in $G_{circ}$.}}}
\label{fig:F3D_2}
\end{figure}

\begin{figure}[tbp]
\begin{center}
\begin{tabular}{cc}
{\includegraphics[scale=0.2,clip=]{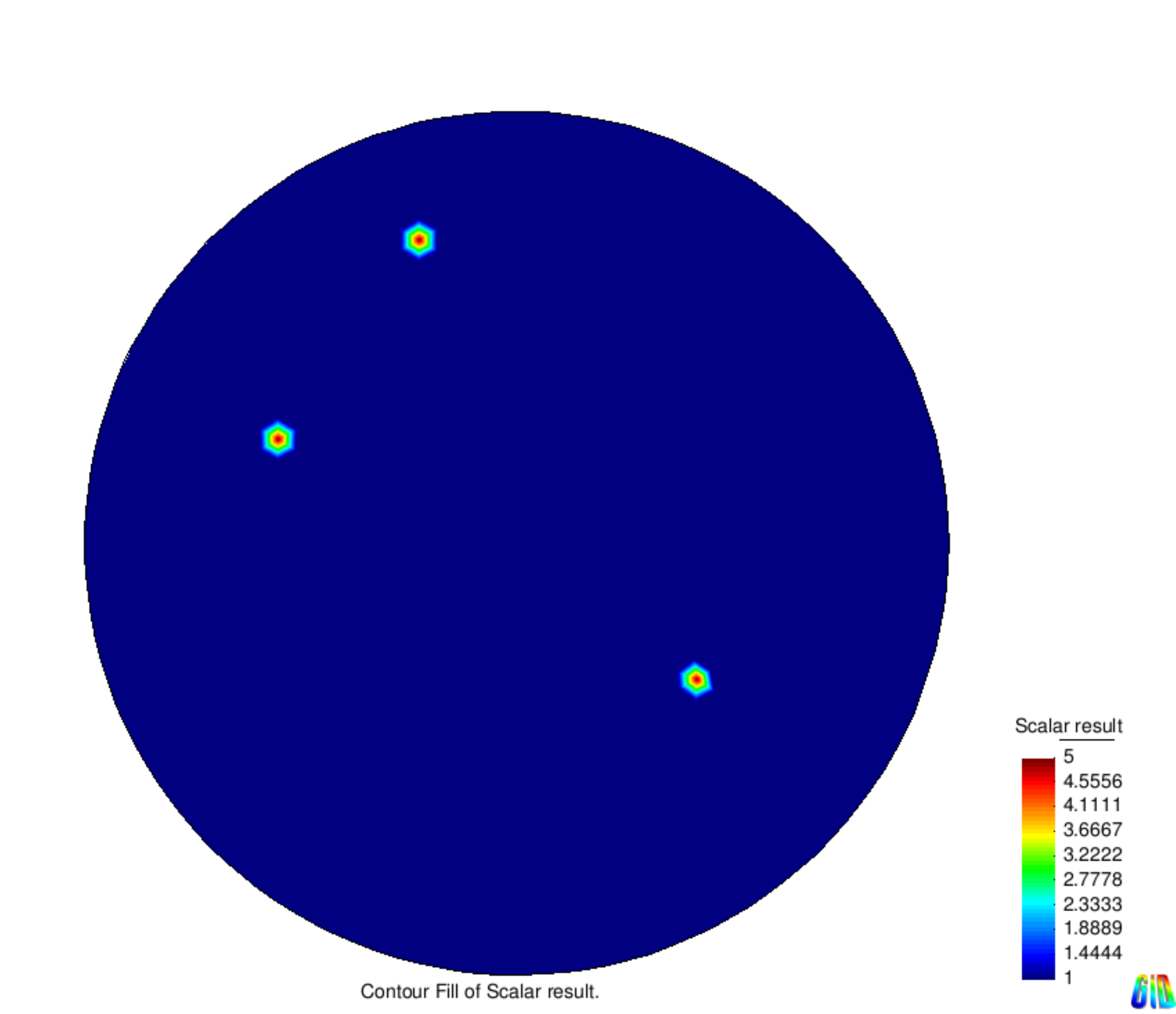}} & {%
\includegraphics[scale=0.2,clip=]{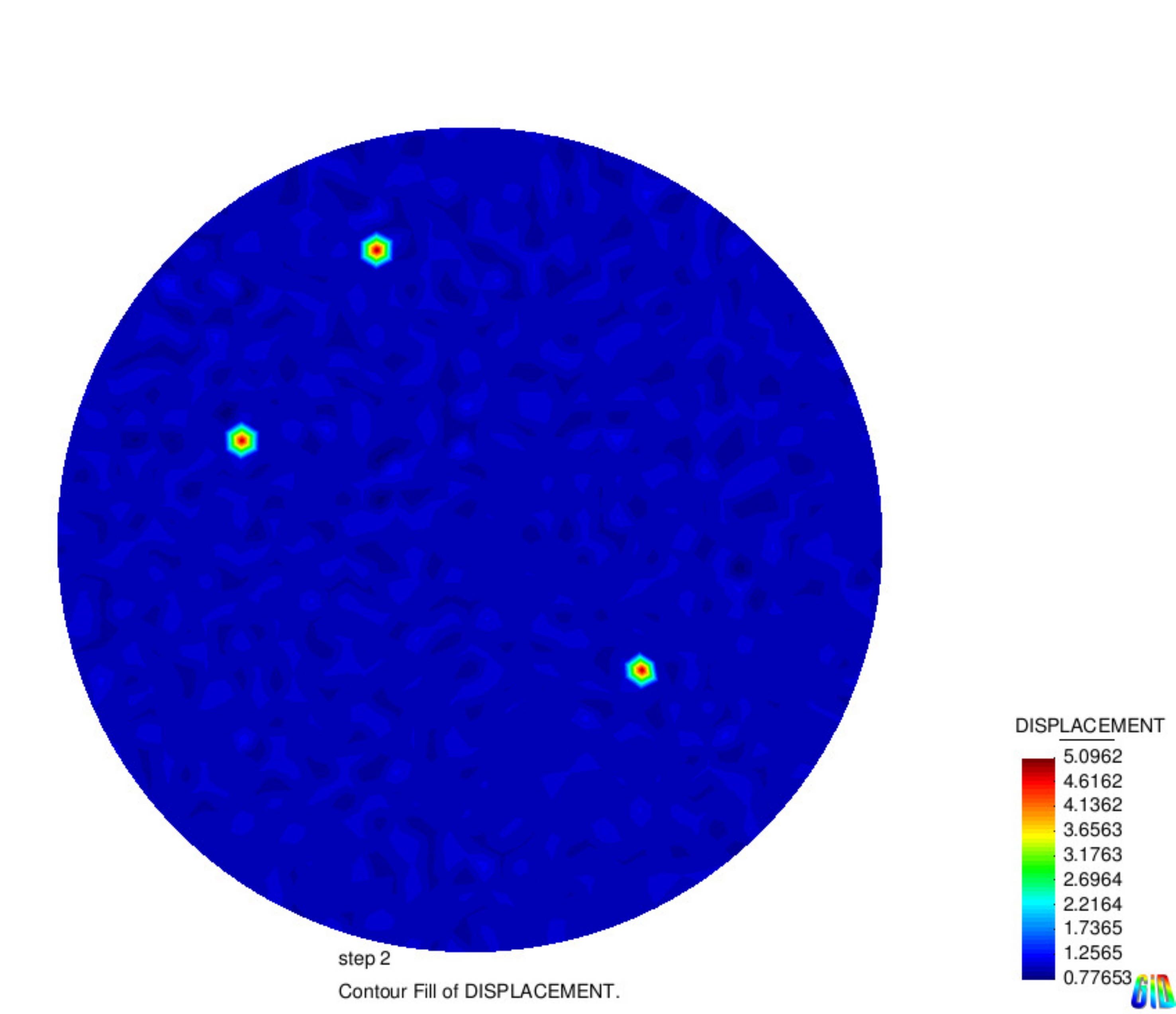}} \\ 
a) exact function $a(x)$ & b) computed at $s=19$ \\ 
{\includegraphics[scale=0.3,clip=]{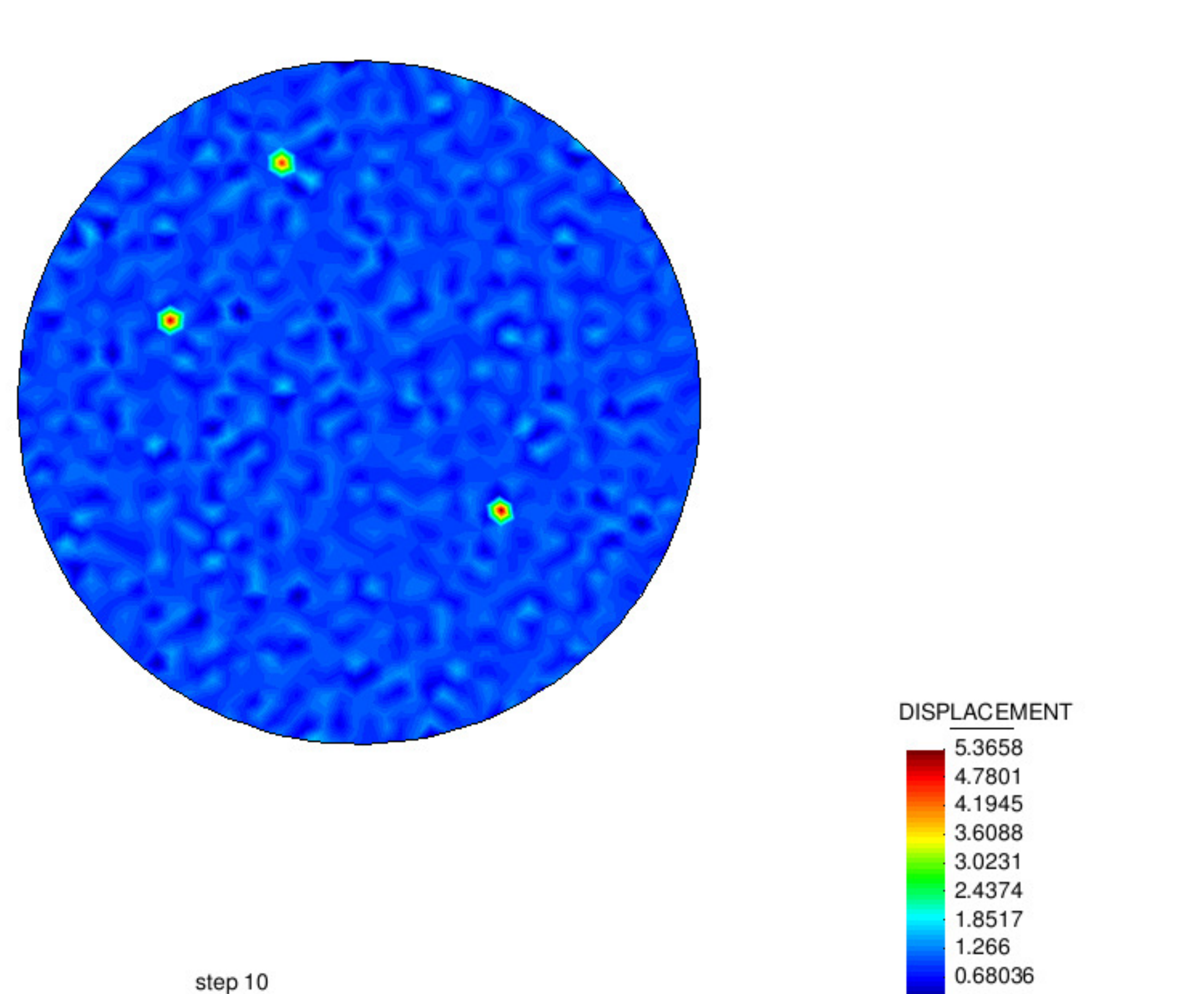}} & {%
\includegraphics[scale=0.3,clip=]{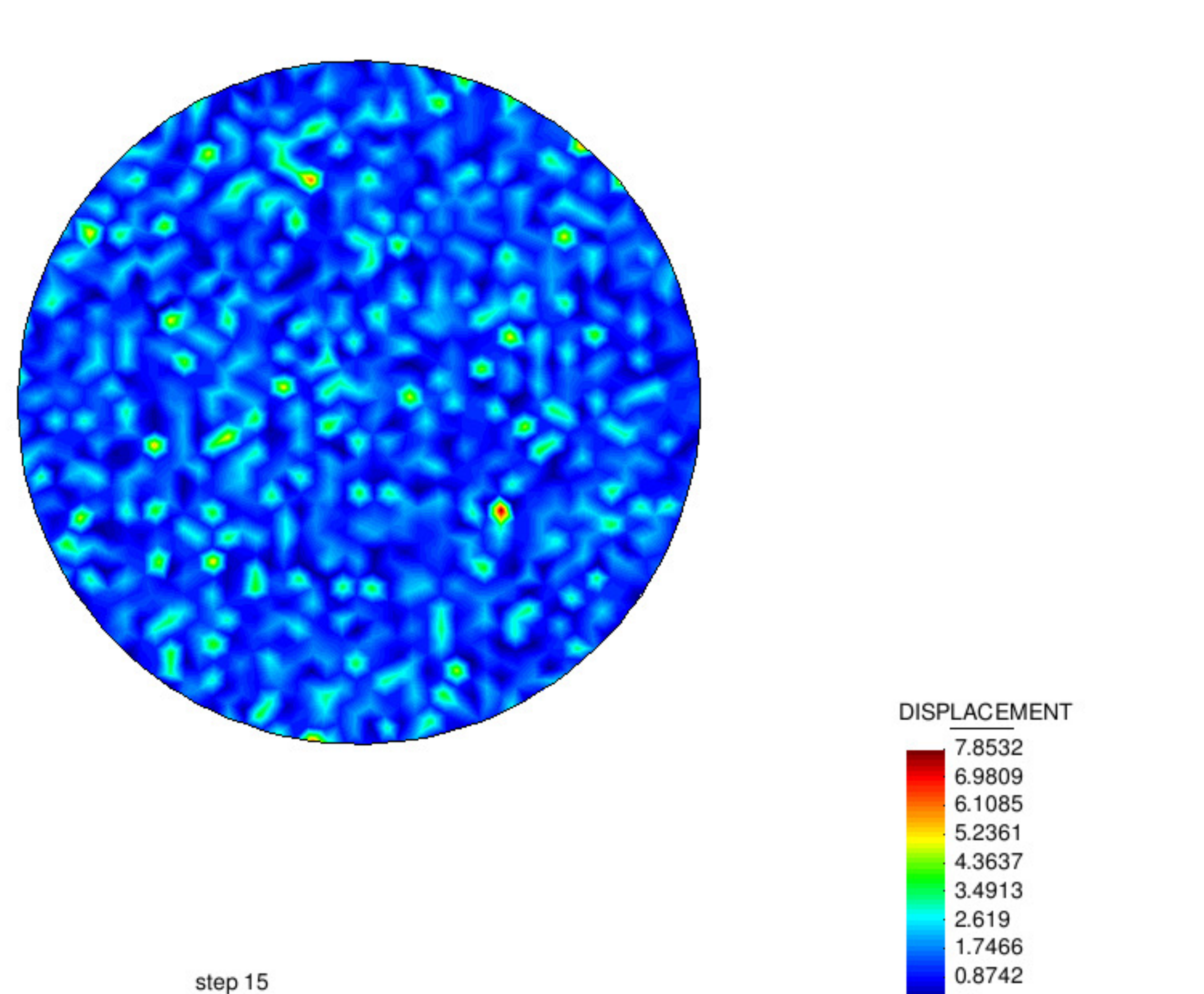}} \\ 
c) computed at $s=10$ &  d) computed at $s=5$
\end{tabular}%
\end{center}
\caption{\emph{\ Test 1: a) The exact location of tumors. b), c), d)
    The reconstructed wave speed function $a(x) $ at different values
    of pseudo frequency $s$ for the the case when the measured
    function $u_{\sigma}(x,t)$ is known inside the domain of interest.
    On b) maximal reconstructed values of this function are $5.09$ in
    three small tumor-like targets. The reconstructed $a(x) =1$
    outside of imaged targets what corresponds to the background
    medium.  Reconstruction presented on b) is highly accurate:
    compare with figure a) where values of the exact function $a(x)$
    inside tumor-like inclusions are $a(x)=5$. However, on d) we observe that
    at pseudo frequency $s=5$ the image is deteriorated.  }}
\label{fig:F5}
\end{figure}

\begin{figure}[tbp]
\begin{center}
\begin{tabular}{cc}
{\includegraphics[scale=0.3,clip=]{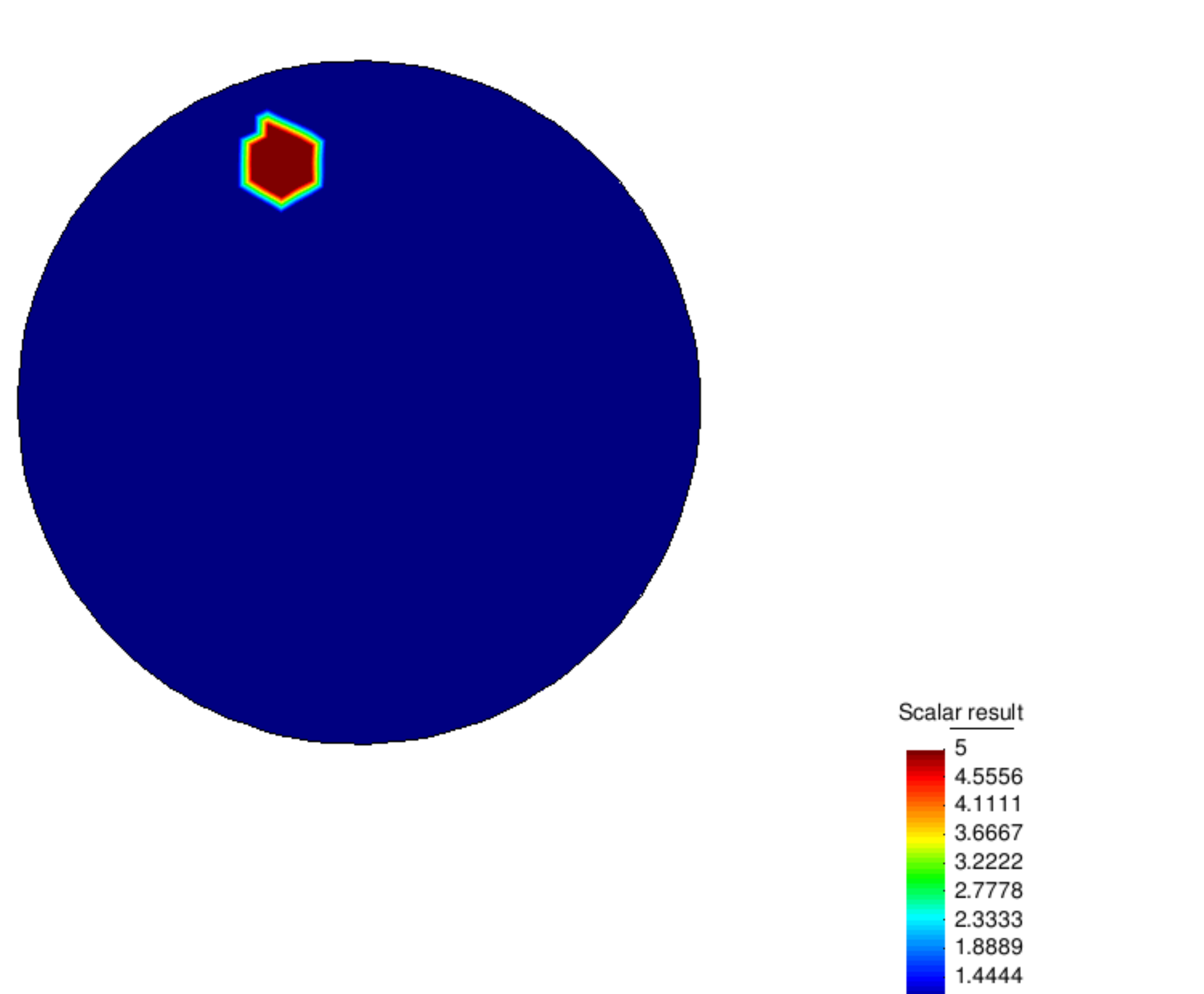}} & 
{\includegraphics[scale=0.3,clip=]{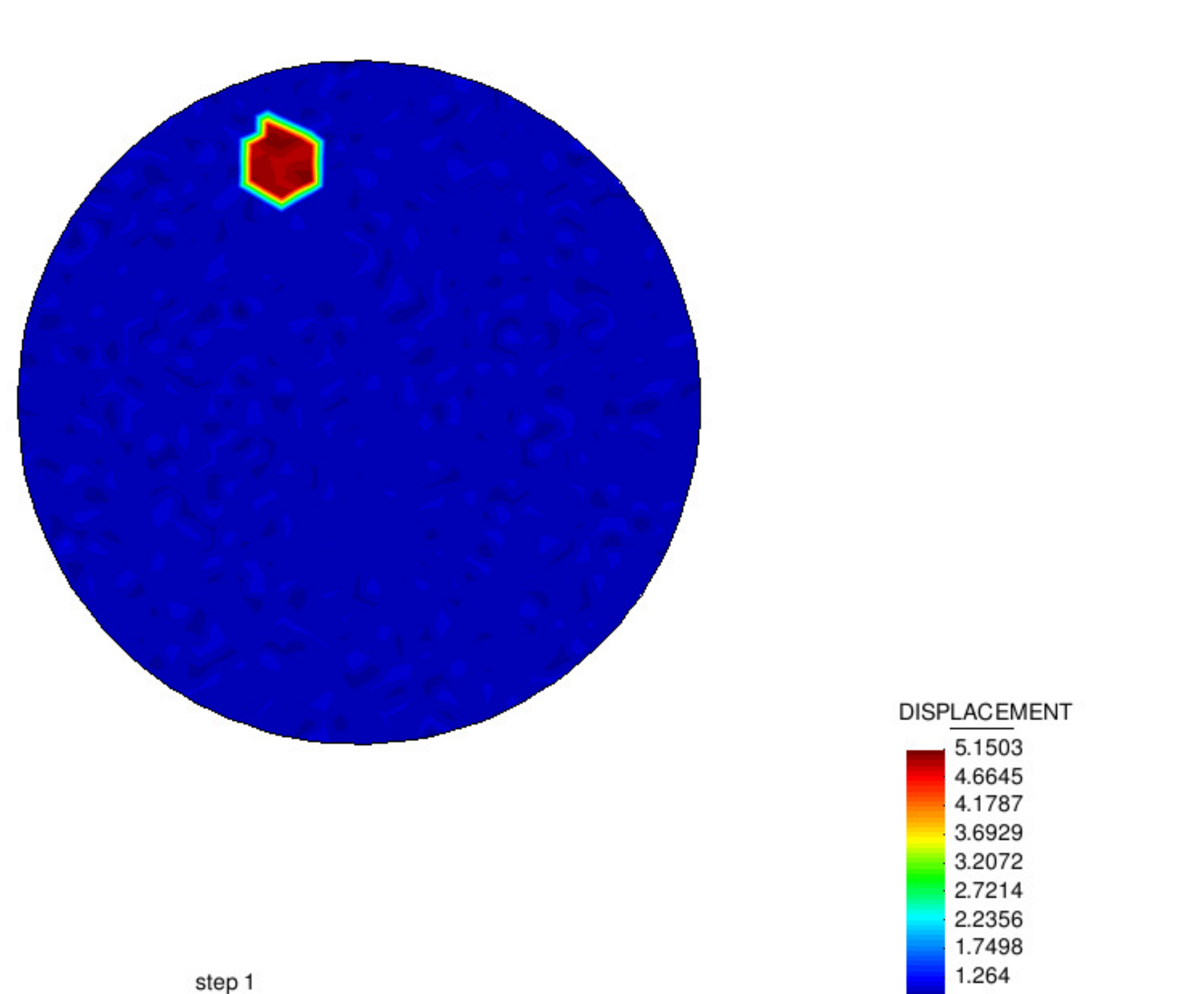}} \\
a) exact function $a(x)$ & b)  computed at $s=19$  \\ 
{\includegraphics[scale=0.3,clip=]{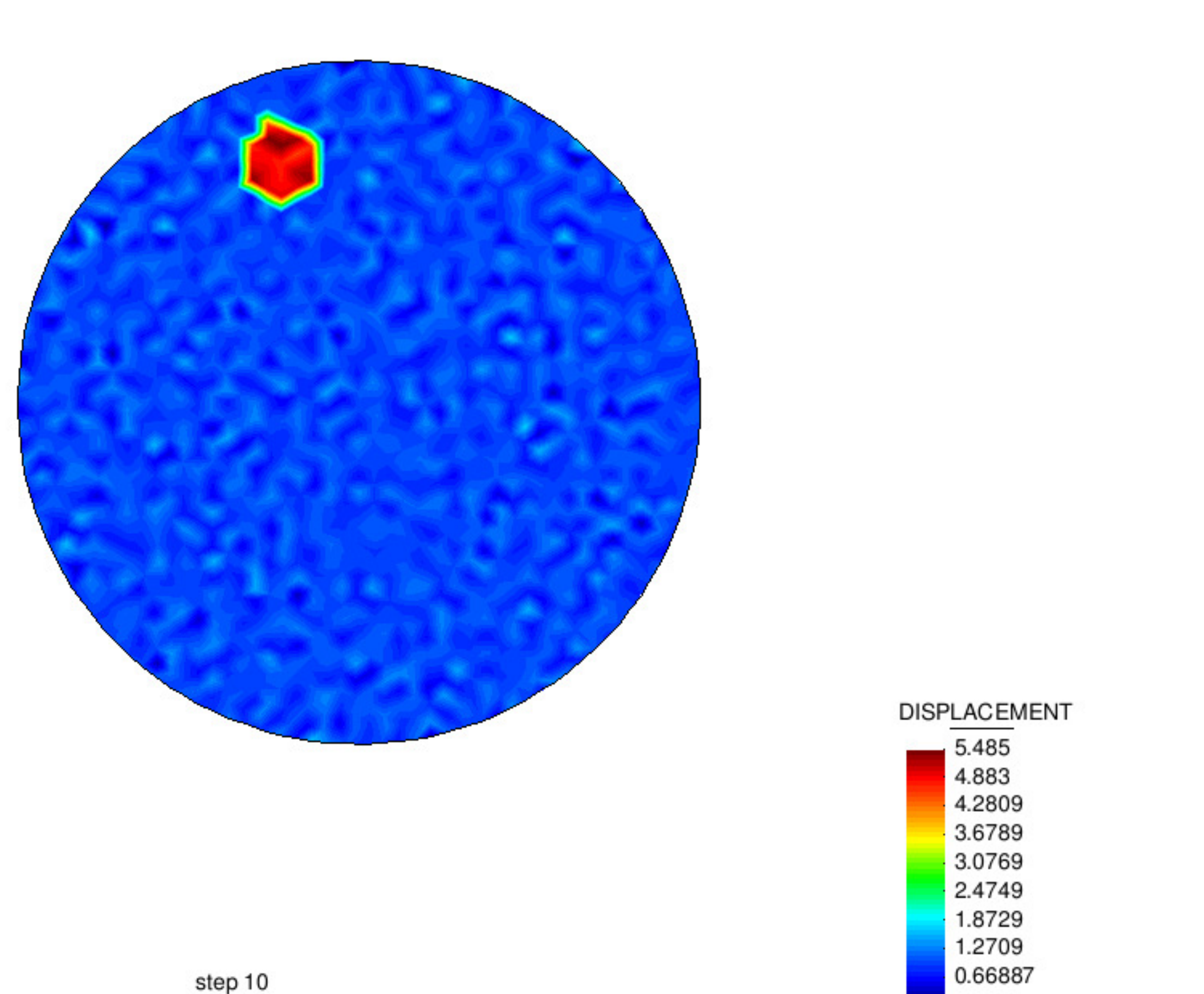}} &
{\includegraphics[scale=0.3,clip=]{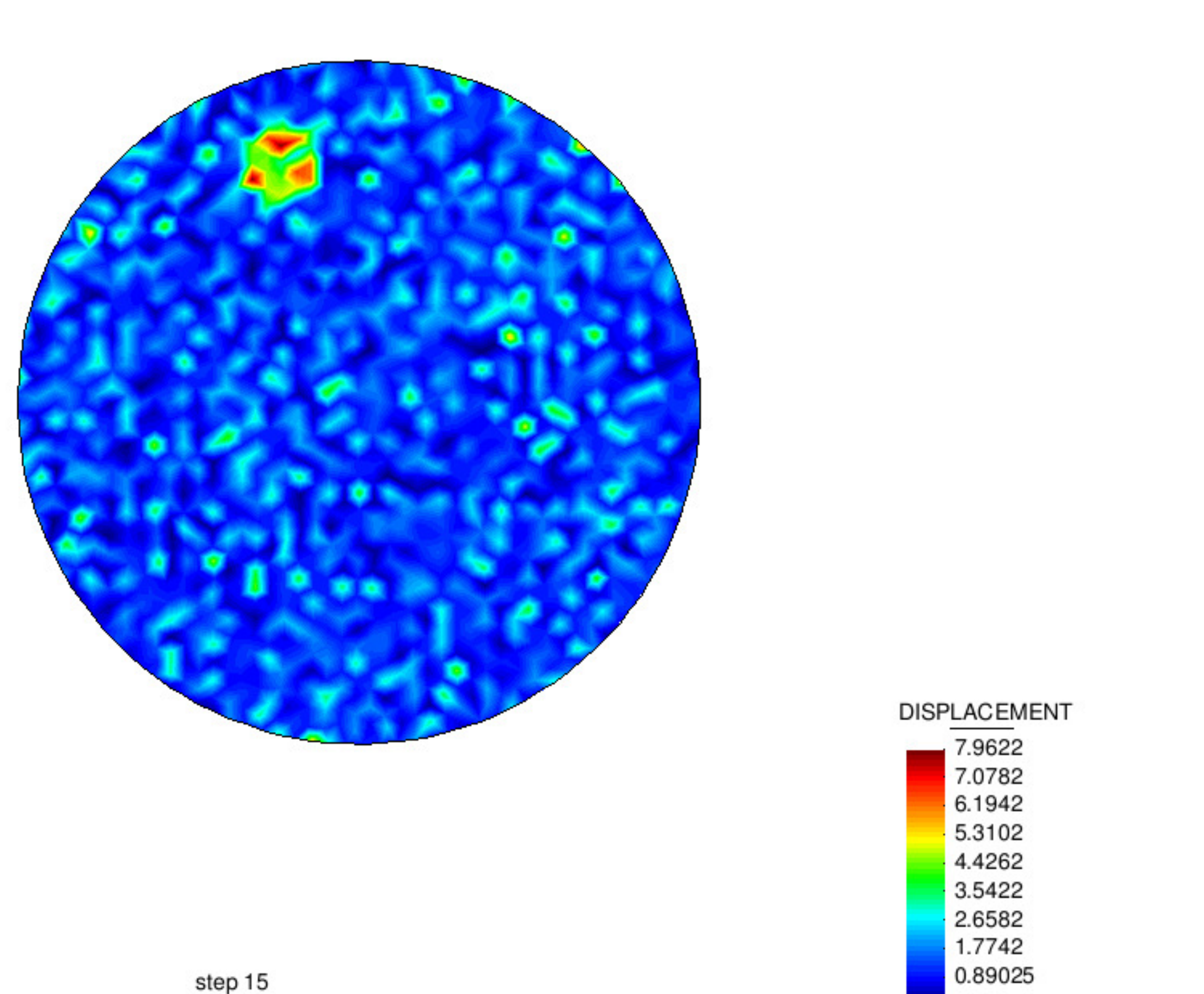}} \\
c)  computed at $s=10$   & d)  computed at $s=5$  \\ 
\end{tabular}
\end{center}
\caption{\emph{\ Test 2: a) The exact location of tumor. b), c), d)
    The reconstructed wave speed function $a(x) $ at different values
    of pseudo frequency $s$ for the the case when the measured
    function $u_{\sigma}(x,t)$ is known inside the domain of interest.
    On b) reconstructed maximal values of this function are  $5.15$ in tumor-like
    target and $a(x) =1$ outside of imaged target what corresponds to
    the background medium. The image presented on b) is highly
    accurate: compare with  figure on a) where the function $a(x)$ in the
    exact tumor-like target has   value $5$. However, on d) we observe that  the image is deteriorated at pseudo frequency $s=5$.  }}
\label{fig:F6}
\end{figure}

\begin{figure}[tbp]
\begin{center}
\begin{tabular}{cc}
{\includegraphics[scale=0.3,clip=]{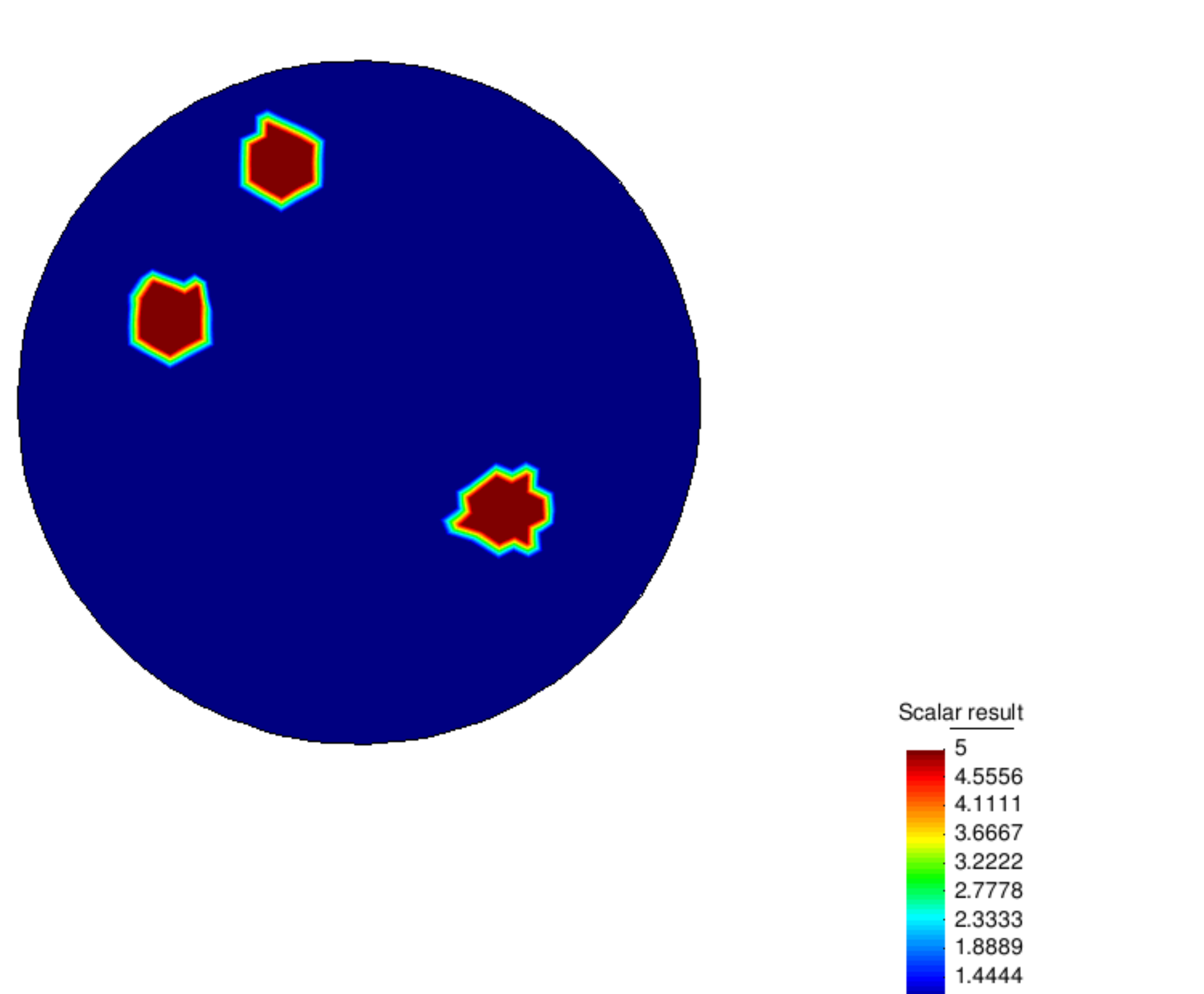}} & 
{\includegraphics[scale=0.3,clip=]{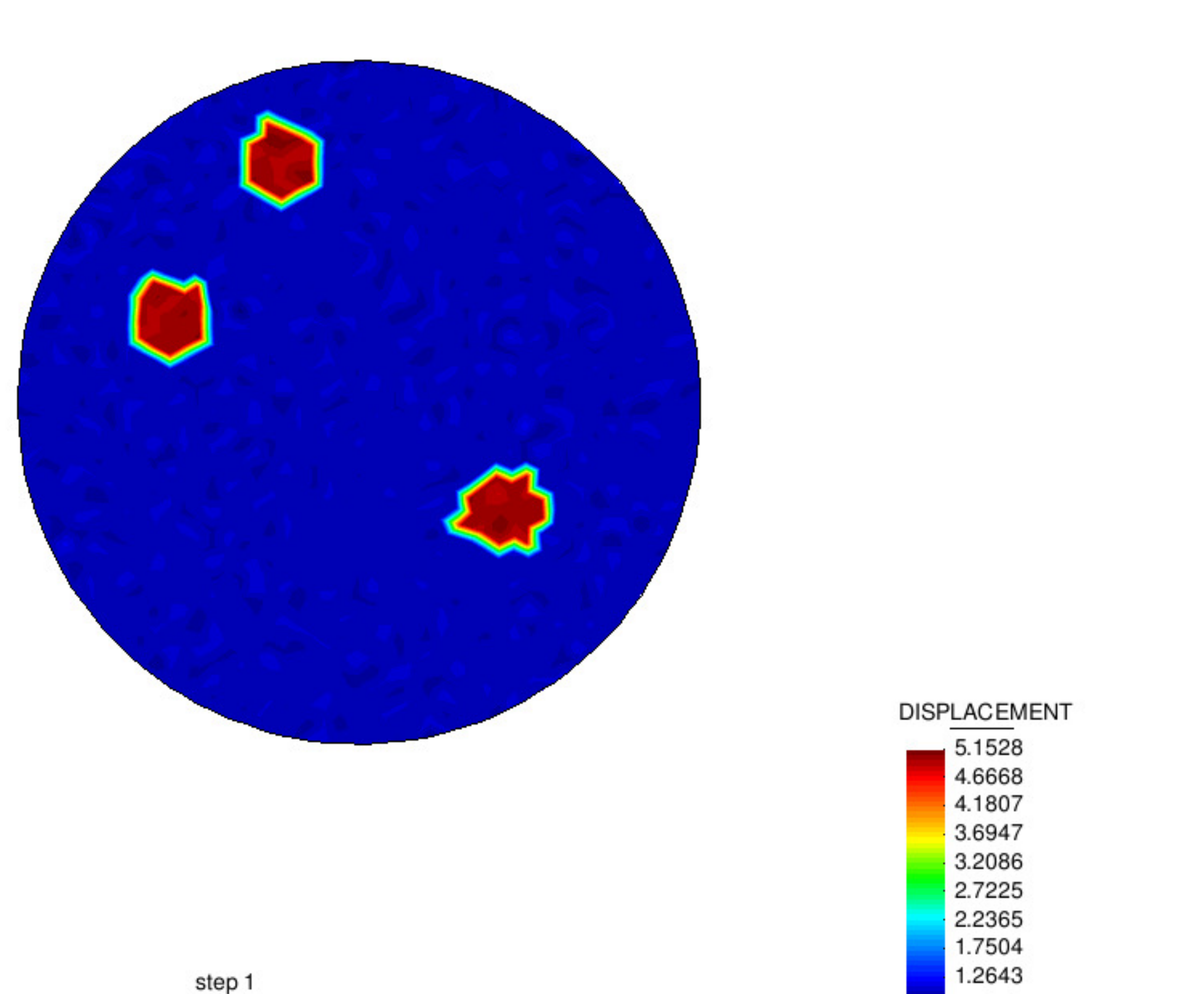}} \\
a) exact function $a(x)$ & b)  computed at $s=19$  \\ 
{\includegraphics[scale=0.3,clip=]{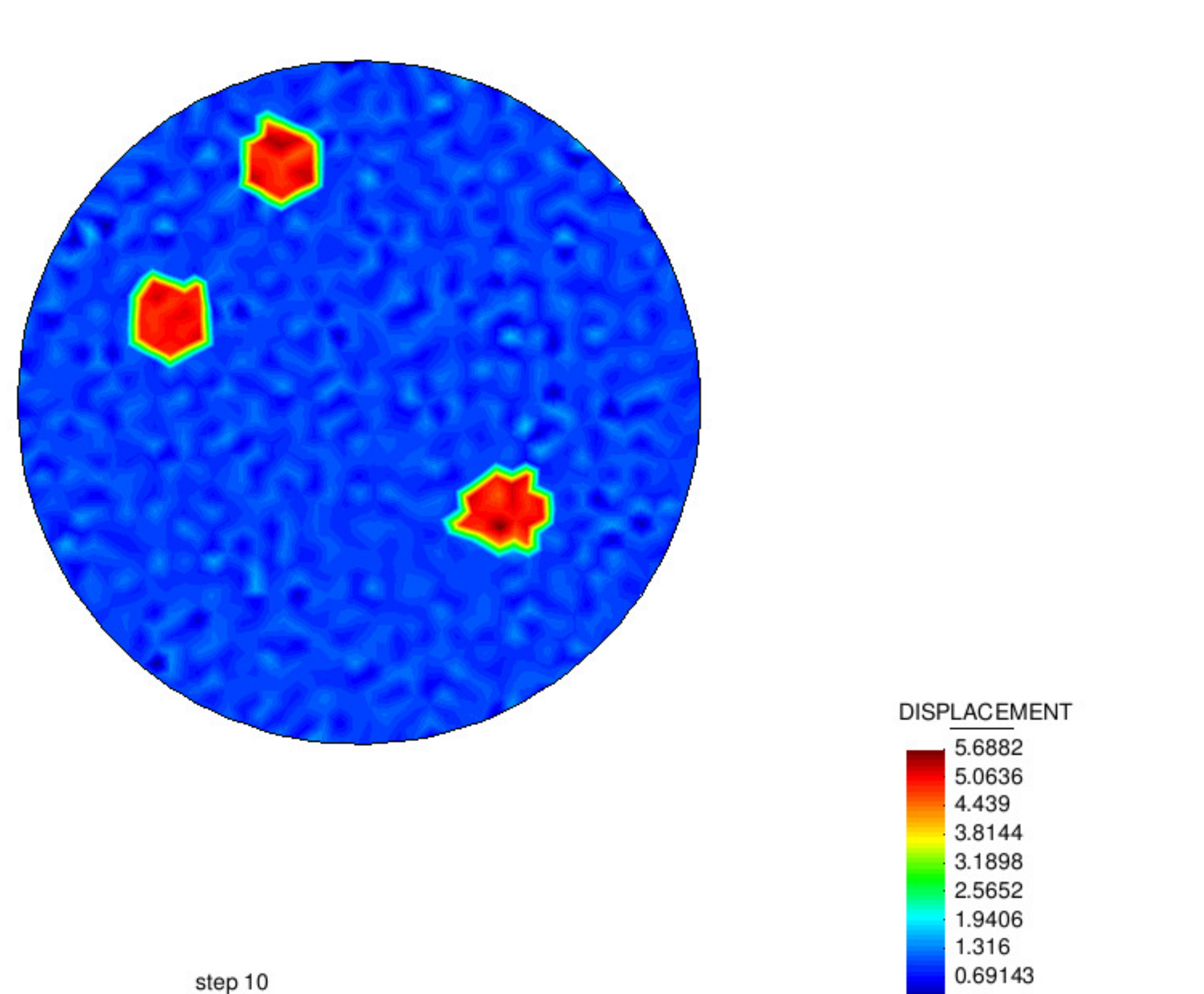}} &
{\includegraphics[scale=0.3,clip=]{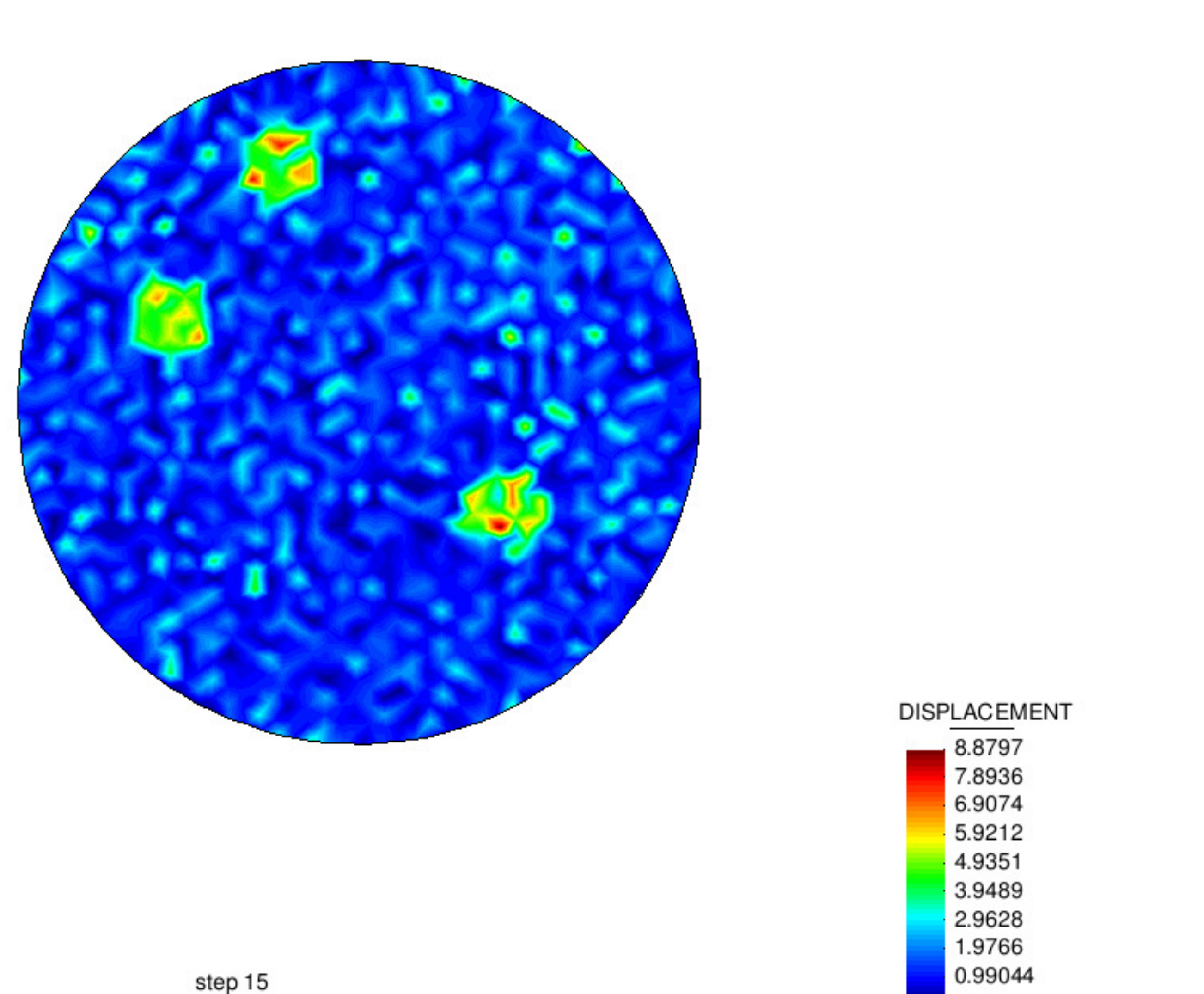}} \\
c)  computed at $s=10$   & d)  computed at $s=5$  \\ 
\end{tabular}
\end{center}
\caption{\emph{\ Test 3: a) The exact location of tumors. b), c), d)
    The reconstructed wave speed function $a(x) $ at different values
    of pseudo frequency $s$.  On b)
    maximal reconstructed values of this function are $5.15$ in tumor-like targets
    and $a(x) =1$ outside of imaged targets what corresponds to the
    background medium. The image is highly accurate: compare with exact image on a) where maximal
    values of the exact function are $5$.  Again, on d) we observe that  the image is deteriorated at pseudo frequency $s=5$. }}
\label{fig:F7}
\end{figure}

\section{Numerical experiments in 2D}

\label{sec:2}

In this section we present the reconstruction of wave speed function
$a(x)$ at different values of pseudo frequency $s$ for the the case
when the measured function $u_{\sigma}(x,t)$ is known inside the
domain of interest. Measuring of the field internally is allowed in  some cases of medical imaging: for example, in
medical resonance elastic imaging  \cite{BB}.
% This is an ideal case which is rarely is available
%in reality. 
%However, at the first step of study of this problem we
%will present reconstruction in this ideal case.

\subsection{Data simulation in 2d}

\label{sec:2.2}

For generation of  data  to solve our CIP, we first solve the forward
problem for the wave equation with known value of a wave speed inside
our domain of interest. Let us define by $G$ the computational domain where
we compute the forward problem.

We simulate the data for the inverse problem using the software
package WavES \cite{waves}. To do that we solve the forward problem
via the hybrid finite element/finite difference (FEM/FDM) method of \cite{BSA}. In this method the computational domain
$G$ is split in two subdomains, $G=G_{FDM}\cup G_{FEM}$, see Figure
~\ref{fig:F1} for these subdomains.  We use structured mesh with FDM
in $G_{FDM}$ and non-structured mesh and FEM in $G_{FEM}=\Omega.$ The
computational domain $G_{FEM}=\Omega$ is also decomposed into two
domains $G_{FEM}= G_{circ} \cup (G_{FEM} \diagdown G_{circ})$, where
$G_{circ}$ is the circular FEM domain where we search tumor-like
inclusions.  The boundary of the rectangle $G$ is $\partial G=\partial
G_{1}\cup \partial G_{2}\cup \partial G_{3}.$ Here, $\partial G_{1}$
and $\partial G_{2}$ are respectively top and bottom sides of the
largest rectangle of Figure \ref{fig:F1}, and $\partial G_{3}$ is the
union of left and right sides of this rectangle.  The space mesh in
$\Omega $ consists of triangles and it consists of squares in
$G_{FDM}$, with the mesh size $\tilde{h}=0.02$ in the overlapping
regions.

We generate the data via solution of the
following forward problem
\begin{equation}
\begin{split}
a\left( x\right) u_{tt}-\Delta u& =0~~~\mbox{in}~G\times (0,T), \\
u(x,0)& =0,~u_{t}(x,0)=0,~\mbox{in}~G, \\
\partial _{n}u\big \vert_{\partial \Omega _{1}}& =f\left( t\right)~%
\mbox{on}~\partial G_{1}\times (0,t_{1}], \\
\partial _{n}u\big \vert_{\partial \Omega _{1}}& =\partial _{t}u ~\mbox{on}%
~\partial G_{1}\times (t_{1},T), \\
\partial _{n}u\big \vert_{\partial G_{2}}& =\partial _{t}u ~\mbox{on}%
~\partial G_{2}\times (0,T), \\
\partial _{n}u\big \vert_{\partial \Omega _{3}}& =0 ~\mbox{on}~\partial
G_{3}\times (0,T).
\end{split}
\label{6.4}
\end{equation}

The plane wave $f\left( t\right) $ is given by
\begin{equation}
f\left( t\right) =\left\{ 
\begin{array}{cc}
\sin \omega t & \text{ for }t\in (0,\frac{2\pi }{%
\omega }]=(0, t_1], \\ 
0  & \text{ for }t\in \left( \frac{2\pi }{\omega },T\right)
\end{array}
\right.  \label{6.5}
\end{equation}
and is initialized at the top boundary $\partial G_{1}$ of the
computational domain $G$ of Figure ~\ref{fig:F1}. The plane wave
    propagates downwards into $G$ and is absorbed at the bottom
    boundary $\partial G_{2}$ for all times $t\in (0,T).$ In addition,
    it is also absorbed at the top boundary $\partial G_{1}$ for times
    $t\in (t_1 ,T)$.  We use first-order absorbing boundary conditions
    \cite{EM}. In our tests we took $\omega =20$ and $T=2$ in
    (\ref{6.5}), see  some simulations of the forward problem on Figure \ref{fig:F3D_1}.
When solving the inverse problem, we assume that the coefficient
$a(x)$ is unknown in the circle $G_{circ} \subset G$ and has a known
constant value $a(x)=1$ in $G\diagdown \Omega$ and in  $G_{FEM}
\diagdown G_{circ}$, see Figure \ref{fig:F1}.

The trace $g\left( x,t\right) $ of the solution $u\left( x,t\right) $
of the wave equation is recorded inside the circle
$G_{circ}$  where we want to reconstruct the function $a(x)$. Next, the coefficient $a(x)$ is 
forgotten, and our goal is to reconstruct this
coefficient for $x \in \Omega $ from the data $\psi \left( x,s\right)$
which are obtained after Laplace transform of the data $g(x,t)$.
We  impose $5 \%$ of additive noise to the data
$u(x,t)$ to get the  measured function
$u_{\sigma}(x,t)$:
\begin{equation}\label{noise}
u_{\sigma}(x_i,t_j) = u(x_i,t_j)[1 + \alpha_j(u_{max}(x_i,t_j) - u_{min}(x_i,t_j))\sigma].
\end{equation}
Here, $u(x_i,t_j)$ is the solution of the problem (\ref{6.4}) at the
mesh point $x_i$ and time moment $t_j \in (0,T)$, $\alpha_j$ is a
random number on the interval $[-1,1]$, $u_{max}(x_i,t_j)$ and
$u_{min}(x_i,t_j)$ are maximal and minimal values of the computed
solution $u(x_i,t_j)$, respectively, and  $\sigma=0.05$ is the level of
the noise.

\subsection{Test 1}

\label{sec:2.1}

We model the problem of imaging of three point-like
tumor inclusions of Figure \ref{fig:F5}-a) as an CIP for the scalar wave equation.
We set  the dimensionless computational domain $G$  as
\begin{equation}
 G =\left( -0.7,0.7\right) \times \left( -0.7, 0.7\right)  \label{G}
\end{equation}
%We
%set
%\begin{equation*}
%\Omega =\left\{ \left( x,y\right) \in \left( -0.52,0.52\right) \text{ m}%
%\times (-0.52,0.52)\text{ m}\right\} ,
%\end{equation*}%
%where \textquotedblleft m" stands for meter. Introducing dimensionless
%spatial variables $\left( x^{\prime },y^{\prime }\right) =\left( x,y\right)
%/\left( 0.1\text{m}\right) $ we
and  the dimensionless domain $G_{FEM}=\Omega$  as
\begin{equation}
\Omega =\left( -0.52,0.52\right) \times \left( -0.52, 0.52\right).  \label{6.2}
\end{equation}
Our domain of interest $G_{circ} \subset G$ where we solve our CIP and
search for tumors, has the center at the point with coordinates
$(0,0)$ and the radius $r=0.4$.  We model our three point-like tumors
$(p_1, p_2, p_3)$ to be located at points of the domain $G_{circ}$
with coordinates
\begin{equation}\label{3points}
\begin{split}
p_1(x_1, y_1):  x_1  &=  -0.090234,  y_1 =  0.280903, \\
p_2(x_2, y_2):  x_2  &= -0.221014,  y_2 = 0.096346, \\
p_3(x_3, y_3): x_3  &= 0.166988,  y_3 = -0.126124.
\end{split}
\end{equation}

Medical experiments show that the relation of the function $a(x)$ in
cancerous tumors to the healthy tissue is $\approx 5$. Thus, we
consider the following relative values of the function $a(x)$ in our
tests
\begin{equation}
 a(x) =\left\{ 
\begin{array}{cc}
1 & \text{ healthy tissue }, \\ 
5  & \text{cancerous tumors}.
\end{array}
\right.  \label{6.3}
\end{equation}

 In Figure \ref{fig:F5}-b) we present reconstruction of three
 tumor-like inclusions of Figure \ref{fig:F5}-a). We use
globally convergent algorithm of section \ref{sec:2.5}
  to get
 reconstructed function $a(x)$ of Figures \ref{fig:F5}-b), c),
 d). Discrete values $a_{n,i}$ at every point $i$ of the computational
 domain $G_{FEM}$ are obtained using formula (\ref{3.109}).
We took pseudo frequency interval $s=[1,19]$ and divided it into subintervals with the step size $\delta s = 1$ for every interval.

Using Figures \ref{fig:F5}-b), c) we observe that we get almost
perfect reconstruction when pseudo frequency $s$ is taken on the
interval $s=[8;19]$. However, for pseudo frequencies on the interval
$s=[1;7]$ we obtain reconstructed function $a(x)$ similar to the one
obtained on Figure \ref{fig:F5}-d). We observe that the image of
Figure \ref{fig:F5}-d) is deteriorated for this value of
pseudo frequency.

\subsection{Test 2}

\label{sec:2.2}

This is the same test as the Test 1 of section \ref{sec:2.1}, only
the goal is image one big tumor-like inclusion of Figure
\ref{fig:F6}-a).

Results are very similar to results of Test 1.
On Figures \ref{fig:F6}-b), c) we observe almost perfect
reconstruction when pseudo frequency $s$ is taken as $s=10$ and
$s=19$. Our numerical tests show that on the interval of
pseudo frequencies $s=[8;19]$ we get reconstruction similar to the
exact one of figure \ref{fig:F6}-b). However, for pseudo frequencies on
the interval $s=[1;7]$ we obtain reconstructed function $a(x)$ similar
to the one obtained on Figure \ref{fig:F6}-d). We observe that the
image of Figure \ref{fig:F6}-d) is deteriorated for this value of
pseudo frequency.

\subsection{Test 3}

\label{sec:2.3}

This is the same test as the Tests 1 and 2 above, only the goal is
image 3 big tumor-like inclusions of Figure \ref{fig:F7}-a). Results
of reconstruction are similar to results of Tests 1 and Test 2 and are
presented on Figure \ref{fig:F7}-b), c), d).

\section{Summary}

\label{sec:9}

 We have applied a finite element method inside the approximately
 globally convergent method of \cite{BK} for explicit reconstruction
 of the coefficient in the hyperbolic equation.  In our numerical
 tests we have used the measured function which was known inside the
 domain of interest. This is possible, for example, in the case of
 magnetic resonance elastography (MRE) which allows measure field
 internally \cite{BB, PO2}. In this work we considered the simplified
 model problem described by the acoustic wave equation instead of the
 elastic one. The elastodynamics system is planned to be considered in
 our future research.
 Results of our numerical examples show
 quantitative and accurate reconstruction of  small tumor-like
 inclusions.

\begin{center}
\textbf{Acknowledgments}
\end{center}

This research was supported by the Swedish
Research Council.


\begin{thebibliography}{99}
 \raggedright



\bibitem{AB}  M.~Asadzadeh and L.~Beilina, A posteriori error analysis
  in a globally convergent numerical method for a hyperbolic
  coefficient inverse problem,  \emph{Inverse Problems}, 26, 115007,
  2010.



\bibitem{BB} P. E. Barbone and J. C. Bamber, Quantitative elasticity
  imaging: what can and cannot be inferred from strain images, \emph{Phys.Med.Biol.}, 47, pp.2147-2164, 2002.


\bibitem{BSA} L. Beilina, K. Samuelsson and K. Åhlander, Efficiency of a
hybrid method for the wave equation. In \emph{\ International Conference on
Finite Element Methods}, Gakuto International Series Mathematical Sciences
and Applications, Gakkotosho CO., LTD, 2001.

\bibitem{BClason} L.~Beilina and C.~Clason, An adaptive hybrid FEM/FDM
method for an inverse scattering problem in scanning acoustic
microscopy, \emph{ SIAM Sci.Comp.}, V.28, I.1, pp.382--402, 2006.

\bibitem{BK} L. Beilina and M.V. Klibanov, \emph{Approximate Global
Convergence and Adaptivity for Coefficient Inverse Problems}, Springer, New
York, 2012.

\bibitem{BK2} L. Beilina and M.V. Klibanov,    A new approximate mathematical model for global convergence for a coefficient inverse problem with backscattering data, J. Inverse and Ill-Posed Problems, 20, pp.513--565, 2012. 


\bibitem{BTKF} L. Beilina, Nguyen Trung Th\`anh, M. V. Klibanov, M. A. Fiddy, Reconstruction from blind  experimental data for an inverse problem for a hyperbolic equation,  \emph{Inverse Problems} 30, 025002, doi:10.1088/0266-5611/30/2/025002, 2014.


\bibitem{BTKJ} L. Beilina, Nguyen Trung Th\`anh, M. V. Klibanov,
  J.Bondestam Malmberg, Reconstruction of shapes and refractive
  indices from backscattering experimental data using the adaptivity, \emph{Inverse Problems} 30, 105007, 2014.

\bibitem{EM} B. Engquist and A. Majda, Absorbing boundary conditions for the
numerical simulation of waves \emph{\ Math. Comp.} 31, 629--651, 1977.

\bibitem{NBKF} Nguyen Trung Th\`{a}nh, L.~Beilina, M.~V. Klibanov and M.~A.
Fiddy, Reconstruction of the refractive index from experimental
backscattering data using a globally convergent inverse method, \emph{SIAM
J. Scientific Computing}, 36 (3), pp.273--293, 2014.

\bibitem{PO2} Wall, David J.N.; Olsson, Peter; van Houten, Elijah E. W., On an
  inverse problem from magnetic resonance elastic imaging, \emph{SIAM
    Journal on Applied Mathematics}, 2011.

\bibitem{waves} WavES, the software package, http://www.waves24.com
\end{thebibliography}
\end{document}